\pdfoutput=1
\RequirePackage{ifpdf}
\ifpdf 
\documentclass[pdftex]{sigma}
\else
\documentclass{sigma}
\fi

\usepackage[all]{xy}
\usepackage{bm}
\usepackage{stmaryrd}
\usepackage{tikz-cd}
\usepackage{tikz}
\usetikzlibrary{arrows}
\usepackage{bbm}

\newtheorem{Theorem}{Theorem}[section]
\newtheorem{Corollary}[Theorem]{Corollary}
\newtheorem{Lemma}[Theorem]{Lemma}
\newtheorem{Proposition}[Theorem]{Proposition}
 { \theoremstyle{definition}
\newtheorem{Definition}[Theorem]{Definition}

\newtheorem{Example}[Theorem]{Example}
\newtheorem{Remark}[Theorem]{Remark} }

\numberwithin{equation}{section}

\usepackage{todonotes}

\newcommand{\ol}{\overline}

\newcommand{\cc}[1]{\mathcal{#1}} 

\newcommand{\CC}{\mathbb{C}}
\newcommand{\ZZ}{\mathbb{Z}}

\newcommand{\PP}{\mathbb{P}}

\newcommand{\UU}{\mathbb{U}}

\newcommand{\RR}{\mathbb{R}}

\newcommand{\bt}{{\bm{t}}}

\newcommand{\torus}{\mathbb{T}}
\newcommand{\FM}{\mathbb{FM}}

\DeclareMathOperator{\ch}{ch}

\DeclareMathOperator{\tot}{tot}

\DeclareMathOperator{\eff}{Eff}

\DeclareMathOperator{\Hom}{Hom}
\DeclareMathOperator{\Gr}{Gr}

\newcommand{\br}[1]{\left\langle#1\right\rangle} 

\newcommand{\op}[1]{\operatorname{#1}}
\newcommand{\F}{F}
\newcommand{\E}{E}
\newcommand{\MB}{M}
\newcommand{\LB}{L}
\newcommand{\z}{z}
\newcommand{\x}{x}
\newcommand{\y}{y}

\begin{document}

\allowdisplaybreaks

\newcommand{\arXivNumber}{2404.12303}

\renewcommand{\PaperNumber}{008}

\FirstPageHeading

\ShortArticleName{Wall Crossing and the Fourier--Mukai Transform for Grassmann Flops}

\ArticleName{Wall Crossing and the Fourier--Mukai Transform\\ for Grassmann Flops}

\Author{Nathan PRIDDIS~$^{\rm a}$, Mark SHOEMAKER~$^{\rm b}$ and Yaoxiong WEN~$^{\rm c}$}

\AuthorNameForHeading{N.~Priddis, M.~Shoemaker and Y.~Wen}

\Address{$^{\rm a)}$~Department of Mathematics, 275 TMCB, Brigham Young University, Provo, UT 84602, USA}
\EmailD{\href{mailto:priddis@math.byu.edu}{priddis@math.byu.edu}}

\Address{$^{\rm b)}$~Department of Mathematics, Colorado State University,\\
\hphantom{$^{\rm b)}$}~1874 Campus Delivery Fort Collins, CO 80523, USA}
\EmailD{\href{mailto:mark.shoemaker@colostate.edu}{mark.shoemaker@colostate.edu}}

\Address{$^{\rm c)}$~School of Mathematics, Korea Institute for Advanced Study, Seoul 02455, South Korea}
\EmailD{\href{mailto:y.x.wen.math@gmail.com}{y.x.wen.math@gmail.com}}

\ArticleDates{Received April 26, 2024, in final form January 27, 2025; Published online February 06, 2025}

\Abstract{We prove the crepant transformation conjecture for relative Grassmann flops over a smooth base $B$. We show that the $I$-functions of the respective GIT quotients are related by analytic continuation and a symplectic transformation. We verify that the symplectic transformation is compatible with Iritani's integral structure, that is, that it is induced by a Fourier--Mukai transform in $K$-theory.}


\Classification{14N35; 14E16; 53D45}

\section{Introduction}

The connection between birational geometry and quantum cohomology has been a source of great activity and interest in Gromov--Witten theory since the early development of the theory. As first described in Ruan's crepant resolution conjecture in \cite{Ru} and later extended to more general crepant transformations in \cite{CIT, CR, LLW, LR}, it is expected that for a general crepant rational map $f\colon X_- \dashrightarrow X_+$, the quantum cohomology of $X_-$ should be equivalent to that of~$X_+$ in a specific sense.

More specifically, suppose we have rational maps of smooth projective varieties
\begin{equation}\label{e1}
 \begin{tikzcd}[sep=small]
 & \tilde X \ar[dl, swap, "f_-"] \ar[dr, "f_+"] &\\
 X_- \ar[rr, dashed, "f"] &&X_+
\end{tikzcd}\end{equation}
such that $f_-^*(K_{X_-}) = f_+^*\bigl(K_{X_+}\bigr)$. The crepant transformation conjecture predicts that generating functions of genus zero Gromov--Witten invariants of $X_-$ and $X_+$ are equal after analytic continuation.
Versions of the conjecture have now been proven in many specific instances, and for certain classes of crepant transformations (see, e.g., \cite{BG, CIJ, ILLW, LR}).

One can further refine the crepant transformation conjecture to take into account Iritani's integral structure (see \cite{IriInt}) or the rational structure of Katzarkov--Kontsevich--Pantev (see \cite{KKP}). With respect to this structure, the correspondence should be compatible with the Fourier--Mukai transform
\[
\FM = {f_+}_* f_-^*\colon \ K^0(X_-) \to K^0\bigl(X_+\bigr).
\]

In this paper, we prove a crepant transformation conjecture for relative Grassmann flops and show it is compatible with the Fourier--Mukai transform.

\subsection{Grassmann flops}
A rich source of examples of
crepant transformations as in \eqref{e1}
arises from variation of GIT. Namely, suppose a reductive group $G$ acts on a variety $V$. Then we can define two GIT quotients
$X_\pm = [V\sslash_\pm G]$,
where $+/-\colon G \to \CC^*$ are characters lying in adjacent maximal chambers with respect to the wall and chamber structure on $\Hom( G, \CC^*) \otimes_\ZZ \RR$.\footnote{For those unfamiliar with the notation, by $[V\sslash_\pm G]$ we mean the quotient stack $[V^{\rm ss}(\pm)/ G]$.} Under certain numerical conditions, the rational map $X_+ \dashrightarrow X_-$ will be crepant.
In this context, the relation between the Gromov--Witten theory of $X_+$ and $X_-$ is known as \emph{wall crossing}, as the variation of GIT between the two varieties involves crossing a codimension-one wall in the $G$-ample cone.

In the present paper, we focus on a particular variation of GIT known as a Grassmann flop. This is one of the simplest examples of variation of GIT with respect to a quotient by a~non-abelian group. In fact, we work in the setting of relative Grassmann flops, as we describe below.

The setup is as follows: given a smooth projective variety $B$, let $\F = \bigoplus_{i=1}^n \MB_i$ and $\E = \bigoplus_{j=1}^n \LB_j$ be sums of line bundles on $B$. For $r < n$, define the vector bundle $V \to B$ to be the total space~of
\[
\mathcal{H}{\rm om}(\F, \cc O_B \otimes \CC^r) \times \mathcal{H}{\rm om}(\cc O_B \otimes \CC^r, \E).
\]
There is a natural (fiberwise) action of $G = {\rm GL}_r$ on $V$.
Let $X_+$ and $X_-$ denote the GIT quotients $V\sslash_\pm G$ associated to the characters $\det^{+1}$ and $\det^{-1}$, respectively. Then
\begin{align*}
X_+ = \tot\bigl(S_+ \otimes p^*\F^\vee\bigr) \to \Gr(r, \E), \qquad
X_- =\tot(S_- \otimes p^*\E) \to \Gr\bigl(r, \F^\vee\bigr),
\end{align*}
where $p\colon \Gr(r, \E) \to B $ \big(resp.\ $p\colon \Gr\bigl(r, \F^\vee\bigr) \to B$\big) denotes the relative Grassmannian of $r$-planes in $\E$ \big(resp.\ $\F^\vee$\big) and $S_\pm$ denotes the rank-$r$ tautological bundle on the associated Grassmann bundle. We will work equivariantly with respect to the action of $\torus = (\CC^*)^{2n}$ induced by the natural action on $E \oplus F$.

\begin{Example}
We give an example of a relative Grassmann flop realizing a variation of GIT of quiver varieties.
Fix integers $0<r< k < n$ and let $B$ be the Grassmannian $\Gr(k, n).$ Let~${F \to B}$ be the trivial rank $k$ bundle $\cc O_B \otimes \CC^k $
and let $E \to B$ be
the rank $k$ tautological bundle on $B$. In this case, 
\[
X_+ = \tot\bigl(\hat S_{r}^{\oplus k}\bigr),
\]
where $\hat S_{r}$ denotes the rank $r$ tautological bundle on $\op{Fl}(r, k, n).$

On the other hand, $X_-$ is the total space of a vector bundle over $\Gr(r, k) \times \Gr(k, n).$ If we let $\bar S_{r}$ denote the tautological rank $r$ bundle from the first factor and $S_{k}$ the tautological rank $k$ bundle from the second factor, then
\[
X_- = \tot\bigl(\bar S_{r}^\vee\otimes S_{k}\bigr).
\]

This example can be realized as an instance of variation of GIT of quiver varieties. Consider the quiver
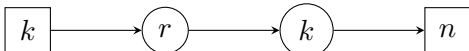
\begin{equation*}
	\begin{tikzpicture}
			\node[draw,
			circle,
			minimum size=0.6cm,
			] (gauge-1) at (0,0){$k$};
		
			\node[draw,
			circle,
 minimum size=0.6cm,
			left=1.2cm of gauge-1
			] (gauge-r) {$r$};

 \node[draw,
 minimum width=0.6cm,
			minimum height=0.6cm,
			left=1.2cm of gauge-r
			] (frame-m) {$k$};

 \node[draw,
			minimum width=0.6cm,
			minimum height=0.6cm,
			right=1.2cm of gauge-1
			] (frame-2) {$n$};

 \draw[-stealth] (frame-m.east) -- (gauge-r.west)
			node[midway,above]{};

			\draw[-stealth] (gauge-r.east) -- (gauge-1.west)
			node[midway,above]{};

 \draw[-stealth] (gauge-1.east) -- (frame-2.west)
			node[midway,below]{};
	\end{tikzpicture}
 \end{equation*}
The associated quiver varieties (for various choices of stability conditions) are quotients of
\[
\op{Hom}\bigl(\CC^{k},\CC^{r}\bigr)\times\op{Hom}\bigl(\CC^{r},\CC^{k}\bigr)\times\op{Hom}\bigl(\CC^{k},\CC^{n}\bigr)
\]
by the natural right action of ${\rm GL}_{r} \times {\rm GL}_{k}$ where ${\rm GL}_{r}$ acts trivially on $\op{Hom}\bigl(\CC^{k},\CC^{n}\bigr)$ and ${\rm GL}_{k}$ acts trivially on $\op{Hom}\bigl(\CC^{k},\CC^{r}\bigr)$. Then $X_+$ and $X_-$ from above are GIT quotients of this representation associated to the characters
\[
+\colon \ (g_1, g_2) \mapsto \det(g_1) \cdot \det(g_2)
\qquad \text{and}\qquad
-\colon \ (g_1, g_2) \mapsto \det(g_1)^{-1} \cdot \det(g_2),
\]
respectively.
\end{Example}

\subsection{Statement of theorem}
With $X_\pm$ as above, consider the infinite-dimensional (symplectic) vector space
\[
\cc H_{X_\pm} = H^*_\torus(X_\pm)\bigl(\!\bigl(z^{-1}\bigr)\!\bigr)[\![Q, q_\pm]\!].
\]
Here $Q$ and $q_\pm$ are Novikov parameters with respect to the base $B$ and fiber directions of $X_\pm \to B$, respectively.

Using the relative quasimap theory of \cite{Oh1}, one can define $I$-functions
$I_{X_\pm}(Q, q_\pm, z)$
lying in $\cc H_{X_\pm}$. These functions are shown to lie on Givental's Lagrangian cone $\cc L_{X_\pm} \subset \cc H_{X_\pm}$ for the respective spaces and may therefore be expressed as generating functions of genus-zero Gromov--Witten invariants of $X_\pm$. See Section~\ref{s:I} for details.

One can then extend $\cc H_{X_\pm}$ by
$\widetilde{\cc H}_{X_\pm} = \cc H_{X_\pm}[\log z]\bigl[z^{-1/2}\bigr]$, and following Iritani in \cite{IriInt}, define the map
\begin{align*}
\Psi\colon \ K^0_\torus(X) \to \widetilde{\cc H}_X, \qquad
A \mapsto z^{- \mu} z^{\rho} \bigl( \hat \Gamma_X \cup (2 \pi {\rm i})^{\op{deg}_0/2} \ch^\torus(A) \bigr),
\end{align*}
where ${\rm i} = \sqrt{-1}$, $\mu$ is the \emph{grading operator}, $\rho = c_1^\torus(X_\pm)$, and $\hat \Gamma_X$ is the Gamma class -- a~characteristic class associated with the Gamma function $\Gamma(1 + z)$.

We prove the following.

\begin{Theorem}[Theorem~\ref{t:st3}]\label{st0}
There exists a linear symplectic isomorphism
$
\UU\colon \widetilde{\cc H}_{X_{-}} \to \widetilde{\cc H}_{X_{+}}
$
such that
the following diagram commutes:
\begin{equation}
 \label{ec1}
\begin{tikzcd}
K^0_\torus(X_{-}) \ar[r, "\FM"] \ar[d, "\Psi_-"] & K^0_\torus(X_{+}) \ar[d, "\Psi_+"] \\
\widetilde{\cc H}_{X_{-}} \ar[r, "\UU"] & \widetilde{\cc H}_{X_{+}}.
\end{tikzcd}
\end{equation}
There is a path $\hat \gamma$ from a neighborhood of $q_+ = 0$ to a neighborhood of $q_- = q_+^{-1} = 0$ such that
\begin{equation}\label{ec2}
 \UU I_{X_{-}} =\widetilde{I_{X_{+}}},
\end{equation}
where $\widetilde{I_{X_{+}}}$ denotes the analytic continuation of $I_{X_{+}}$ along $\hat \gamma.$

\end{Theorem}

Theorem~\ref{st0} generalizes previous known cases of the crepant transformation conjecture in two directions (see Section~\ref{ss:otherwork} for more detailed comparisons with previous work). First, to our knowledge it is one of the first examples
of a wall crossing result based on analytic continuation of generating functions for a variation of GIT quotients by a non-abelian group -- in this case~${\rm GL}_r$. Second, it provides an example of the conjecture for fiber bundles over an arbitrary smooth base~$B$. Thus, in the special case of the quotient by ${\rm GL}_1 = \CC^*$, where the two GIT quotients are toric bundles over B, our result generalizes a particular case of the main result in \cite{CIJ}.

We expect these techniques to be applicable beyond Grassmann flops, to more general variation of GIT quotients with respect to actions of ${\rm GL}_n$. In this way, the current paper may be viewed as a proof of concept. It would also be very interesting to extend these methods to quotients by other non-abelian groups. We hope to revisit this question in future work.

\subsection{Strategy of proof}
The proof relies on the abelian/non-abelian correspondence in Gromov--Witten theory, which allows us to construct the $I$-functions $I_{X_\pm}$ for the GIT quotients $X_\pm$ in terms of $I$-functions of the associated abelian quotients
\[
X_{T, \pm}:= V \sslash_{\pm} T,
\]
where $T$ is the maximal torus in $G$ consisting of diagonal matrices.
More precisely, the function $I_{X_\pm}(Q, q_\pm, z)$ is obtained as the specialization of
\[
I_{X_{T, \pm}}^{e_\torus(\mathfrak g /\mathfrak t)}(Q, q_{1,\pm}, \dots, q_{r,\pm}, z) \qquad \text{at $q_{1,\pm}= \cdots = q_{r,\pm} = q_\pm$},
\]
where \smash{$I_{X_{T, \pm}}^{e_\torus(\mathfrak g /\mathfrak t)}$} is a twist of the $I$-function for the associated abelian quotient by the root bundle for $G$, which we denote by $\mathfrak g /\mathfrak t$. See Section~\ref{s:I} for details.

The technical heart of the proof relies on computing the analytic continuation of $I_{X_{+}}$ along the path $\hat \gamma$ and using this to determine $\UU$. It is not clear how to compute $\widetilde{I_{X_{+}}}$ directly. Instead, we~compute the analytic continuation of
\smash{$I_{X_{T, +}}^{e_\torus(\mathfrak g /\mathfrak t)}(Q, q_{1,+}, \dots, q_{r,+}, z)$} along each of the variables~$q_{k,+}$ independently. That is to say, we analytically continue \smash{$I_{X_{T, +}}^{e_\torus(\mathfrak g /\mathfrak t)}(Q, q_{1,+}, \dots, q_{r,+}, z)$} along a~concatenated path $\gamma_1 \star \cdots \star \gamma_r$, where each $\gamma_k$ leaves the variables $q_{i, +}$ constant for $i \neq k$, and goes from $0$ to $\infty$ in the $q_{k, +}$ direction. We then \emph{homotope} the path $\gamma_1 \star \cdots \star \gamma_r$ to a~path~$\hat \gamma$ which lies entirely in the locus $q_{1,+}= \cdots = q_{r,+},$ thereby obtaining $\widetilde{I_{X_{+}}}$.

Finally, we must verify that the linear map $\UU$ identifying $I_{X_{-}}$ with $\widetilde{I_{X_{+}}}$ is compatible with the Fourier--Mukai transform. This involves an explicit computation of $\FM$ on generators of~$K^0_\torus(X_-)$. The proof that the map $\UU$ in \eqref{ec1} satisfies \eqref{ec2} then involves a nontrivial identity of anti-symmetric functions (Lemma~\ref{l:sym}).

\subsection{Relation to other work}\label{ss:otherwork}

The structure of the proof follows that given in \cite{ChIR,ChR, CIJ, CIT}. In particular, the formulation of the correspondence and the proof via localization is closely related to~\cite{CIJ}. The new ingredients in the current article are the homotopy of the path of analytic continuation, the generalization to the relative setting, and the combinatorial identity verifying agreement with $\FM$, which does not appear in the case of abelian variation of GIT.

In the special case of the quotient by ${\rm GL}_1 = \CC^*$, where the two GIT quotients are toric bundles over $B$, (a compactification of) this example was studied as the \emph{local model} in \cite{LLW2} and a similar result was obtained at the level of $D$-modules and quantum cohomology using different methods. It is possible that the comparison of $I$-functions in this paper may be obtained from their results, although it is not immediately apparent how to pursue such a strategy. On the other hand, our comparison with the $K$-theoretic Fourier--Mukai transform appears to be new even in this case.

The Gromov--Witten theory of spaces related by variation of GIT quotients by nonabelian groups has also been studied in \cite{GW} in the very general setting of GIT quotients of a projective variety by a complex reductive group. The setting of the current paper is different, as our GIT quotients are of quasi-projective varieties, and the form of the comparison is also different.
The generating functions being compared in loc.~cit.\ are different than in this paper, and the difference between the two generating functions is expressed as an explicit sum of residues rather than via analytic continuation.

In the final stages of writing this paper, we learned of a related result by Lutz--Shafi--Webb (see \cite{LSW}).
While both papers use the idea of homotoping the path of analytic continuation, there are significant differences in the two results and their proofs.
The most important such difference is in the analysis of the symplectic transformation~$\UU$. In this paper, we compare the map~$\UU$ to the Fourier--Mukai transform $\FM$. In \cite{LSW}, they instead compare~$\UU$ with the symplectic transformation of the associated abelian wall crossing. This leads to two different proofs that the map~$\UU$ is a symplectic isomorphism and has a well-defined non-equivariant limit. See Remark~\ref{r:bigcomp} for an explicit comparison of approaches.
The two papers are in fact complementary, in that our respective results may be combined to show that the abelian/non-abelian correspondence is compatible with the Fourier--Mukai transform in $K$-theory. Such a compatibility is not \textit{a priori} obvious, as it is not clear how best to relate the kernel $\tilde X$ of the Fourier--Mukai transform of the Grassmann flop to the kernel of Fourier--Mukai transform of the associated abelian wall crossing.
We hope to revisit this issue in the future to obtain a more direct and conceptual understanding of this compatibility.

\section{Geometric setup}\label{gs}
Let $\E = \bigoplus_{j=1}^n \LB_j$ and $\F = \bigoplus_{i=1}^n \MB_i$ be sums of line bundles over a smooth projective variety $B$. Choose $r < n$, and define the vector bundle $V \to B$ to be the total space of
\[
\mathcal{H}{\rm om}(\F, \cc O_B \otimes \CC^r) \times \mathcal{H}{\rm om}(\cc O_B \otimes \CC^r, \E).
\]
This vector bundle comes equipped with a (right) action of $G = \op{Aut}(\CC^r) = {\rm GL}_r$ given fiberwise~by
\[
(X, Y)\cdot g = \bigl( g^{-1}X, Yg\bigr).
\]
Let $X_+$ and $X_-$ denote the GIT quotients $V\sslash_\pm G$ associated to the characters $\det^{+1}$ and $\det^{-1}$, respectively. Then
\begin{align*}
X_+ = \tot\bigl(S_+ \otimes p^*\F^\vee\bigr) \to \Gr(r, \E) ,\qquad
X_- =\tot(S_- \otimes p^*\E) \to \Gr\bigl(r, \F^\vee\bigr),
\end{align*}
where $p\colon \Gr(r, \E) \to B $ \big(resp.\ $p\colon \Gr\bigl(r, \F^\vee\bigr) \to B$\big) denotes the relative Grassmannian of $r$-planes in $\E$ \big(resp.\ $\F^\vee$\big) and $S_\pm$ denotes the rank-$r$ tautological bundle on the associated Grassmann bundle. In other words, the fiber of
\[
\Gr(r, \E) \to B
\]
over a closed point $b \in B$ is the Grassmannian $\Gr(r, \E_b)$ where $\E_b$ is the fiber of $\E \to B$ over $b$.

The birational map
\[
f\colon \ X_{-}\dashrightarrow X_{+}
\]
is a relative version of what is known as a Grassmann flop. We construct a resolution of this map in Section~\ref{s:FM} following \cite{BLV}, and compute the associated Fourier--Mukai transform in $K$-theory.

\subsection{Cohomology}\label{s:coh}
Let $(\CC^*)^n$ act on $\F$ by scaling, such that the $i$th factor scales $\MB_i$. Let a different copy of $(\CC^*)^n$ act on $\E$ similarly.
This induces an action of the torus $\torus:= (\CC^*)^{2n}$ on both $X_+$ and $X_-$. Unless otherwise specified, all cohomology groups are $\torus$-equivariant. We define $\z_i = c_1^\torus\bigl(\MB_i^\vee\bigr)$ and $\x_i = c_1^\torus\bigl(\LB_i^\vee\bigr)$.

Let $R = \CC^r$ denote the right representation of $G$ (viewing $\CC^r$ as row vectors), given by
\begin{equation} \label{e:srep}
v\cdot g = vg
.\end{equation}
The $G$-equivariant vector bundle $R \times V \to V$ induces a vector bundle on both $X_+$ and $X_-$ which, by abuse of notation, we also denote by $R$. On $X_-$ it is (the pullback of) the tautological bundle $S_- \to \Gr\bigl(r, \F^\vee\bigr)$ and on $X_+$ it is (the pullback of) the dual of the tautological bundle $S_+ \to \Gr(r, \E)$. Denote the equivariant Chern roots of the dual $R^\vee$ by $\y_1, \dots , \y_r$.

Let
$\pi_{\pm}\colon X_\pm \to B$ denote the projection. To simplify notation, we will denote $\pi_\pm^*(x_i)$ (resp.\ $\pi_\pm^*(z_j)$) simply by $x_i$ (resp.\ $z_j$) when no confusion will result.
The equivariant cohomology of $X_+$ (resp.\ $X_-$), denoted $H^*_\torus(X_{\pm})$, is generated as an algebra by elements of
$H^*(B)[\z_i, \x_i]_{i=1}^n$
and $\sigma_1, \dots, \sigma_r$, the elementary symmetric polynomials in $\y_1, \dots , \y_r$, subject to the relations:
\begin{equation*}
I_+ = \left\{ \left[ \frac{\prod_{i=1}^n(1 - \x_i)}{\prod_{j=1}^r (1 + \y_j)}\right]_{l} = 0\ \middle| \ l> n-r \right\}
\end{equation*}
in $H^*_\torus\bigl(X_+\bigr)$, and
\begin{equation*}
I_- = \left\{ \left[ \frac{\prod_{i=1}^n(1 + \z_i)}{\prod_{j=1}^r (1 - \y_j)}\right]_{l} = 0 \ \middle| \ l> n-r \right\}
\end{equation*}
in $H^*_\torus(X_-)$. (See \cite[Example 14.6.6]{Fulton} for an equivalent presentation in the context of Chow rings, the above presentation follows by the same argument.) In the above formulas, $\left[ - \right]_{l}$ denotes the degree $l$ part. We will denote by $R_\torus$ the ring $H^*_\torus(pt)$, and by $S_\torus$ its localization with respect to non-zero homogeneous elements. We will also make use of the completion
$\hat S_\torus$, defined by
\[
\hat S_\torus = \biggl\{ \sum_{d \in \ZZ} a_d \, \bigg| \, \text{there exists $d_0 \in \ZZ$ such that $a_d = 0$ for $d < d_0$}\biggr\},
\]
where $a_d$ lies in the degree $d$ homogeneous part of $S_\torus.$
We will often work with the localized equivariant cohomology ring $H^*_\torus(X_\pm) \otimes_{R_\torus} S_\torus,$ or its completion $H^*_\torus(X_\pm) \otimes_{R_\torus} \hat S_\torus,$

The connected components of the $\torus$-fixed locus of $X_+$ (resp.\ $X_-$) are indexed by subsets of $\{1, \dots, n\}$ of size $r$, corresponding to choosing an $r \times r$ minor of $\mathcal{H}{\rm om}(\cc O_B \otimes \CC^r, \E)$ (resp.\ $\mathcal{H}{\rm om}(\F, \cc O_B \otimes \CC^r)$). Thus each connected component is isomorphic to $B$. Let $D_{+}$ (resp.\ $D_{-}$) denote the set indexing these fixed loci.
We may represent an element of $D_\pm$ by a function $\delta\colon \{1, \dots, r\} \to \{1, \dots, n\}$ such that
$\delta(i)<\delta(j)$ for all $i<j$.
 We will further denote $\delta(i)$ by $\delta_i$.

\begin{Lemma}\label{l:rest}\quad
\begin{enumerate}\itemsep=0pt\itemsep=0pt
\item[$1.$] For $\delta^- \in D_{-}$, let $B_{\delta^-} $ denote the associated fixed locus, isomorphic to $B$. Then
\[ R|_{B_{\delta^-}} \cong \bigoplus_{i=1}^r \MB_{\delta^-_i}^\vee.\]
\item[$2.$] For $\delta^+ \in D_{+}$, let $B_{\delta^+} $ denote the associated fixed locus, isomorphic to $B$. Then
\[R|_{B_{\delta^+}} \cong \bigoplus_{j=1}^r \LB_{\delta^+_j}^\vee.\]
\end{enumerate}
\end{Lemma}

\begin{proof}
We will prove the first statement. The fixed locus may be expressed as the GIT quotient
\[
\mathcal{H}{\rm om}\biggl(\bigoplus_{i=1}^r \MB_{\delta^-_i}, \cc O_B \otimes \CC^r\biggr) \sslash_{\theta_{-}} G.
\]
There is a tautological map of vector bundles
$q^*\bigl(\bigoplus_{i=1}^r \MB_{\delta^-_i}\bigr) \to R^\vee$
over the stack
\[
\biggl[\mathcal{H}{\rm om}\biggl(\bigoplus_{i=1}^r \MB_{\delta^-_i}, \cc O_B \otimes \CC^r\biggr) / G\biggr],
\]
where \smash{$q\colon \bigl[\mathcal{H}{\rm om}\bigl(\bigoplus_{i=1}^r \MB_{\delta^-_i}, \cc O_B \otimes \CC^r\bigr) / G\bigr] \to B$} is the projection.
The locus where this map is an~isomorphism is exactly the $\theta_{-}$-stable locus of $\mathcal{H}{\rm om}\bigl( \bigoplus_{i=1}^r \MB_{\delta^-_i}, \cc O_B \otimes \CC^r\bigr)$.
\end{proof}

\begin{Corollary}\label{c:rest}\quad
\begin{enumerate}\itemsep=0pt
\item[$1.$] In $H^*_\torus(X_-)$, up to a choice of ordering, the restriction of $\y_i$ to $B_{\delta^{-}}$ is $-\z_{\delta^{-}_i}$.
\item[$2.$] In $H^*_\torus\bigl(X_+\bigr)$, up to a choice of ordering, the restriction of $\y_i$ to $B_{\delta^{+}}$ is $-\x_{\delta^{+}_i}$.
\end{enumerate}
\end{Corollary}

\subsection{The associated abelian quotients}\label{s:aa}
Let $T \subset G$ denote the maximal torus in $G$ consisting of diagonal matrices. The characters $\det^\pm$ restrict to characters of $T$. We define the abelian quotients $X_{T, \pm}$ associated to $X_\pm$ to be the GIT quotients $V \sslash_\pm T$. They are fiber products over $B$ of $r$ copies of $\tot\bigl(\cc O_{\PP(\F^\vee)}(-1) \otimes p^*\E\bigr)$ and $r$ copies of $\tot\bigl(\cc O_{\PP(\E)}(-1) \otimes p^*\F^\vee\bigr)$, respectively:
\begin{align*}
&X_{T, {-}}= \tot \bigl( \cc O_{\PP(\F^\vee)}(-1) \otimes p^*\E\bigr) \times_B \cdots \times_B \tot \bigl( \cc O_{\PP(\F^\vee)}(-1) \otimes p^*\E\bigr),\\
&X_{T, {+}}= \tot \bigl( \cc O_{\PP(\E)}(-1) \otimes p^*\F^\vee\bigr) \times_B \cdots \times_B \tot \bigl( \cc O_{\PP(\E)}(-1) \otimes p^*\F^\vee\bigr).
\end{align*}

There is a series of intermediate GIT quotients defined as follows. For $0\leq j \leq r$, let $\theta_j \in \Hom(T, \CC^*)$ denote the character
\[
\theta_j ((t_1, \dots, t_r)) = \prod_{i=1}^{j} t_i \cdot \prod_{i=j+1}^r t_i^{-1},
\]
and let $X_{T, j}$ denote the GIT quotient $ V \sslash_{\theta_j} T$:
\[
X_{T, j} = \prod_{i=1}^{j} \tot \bigl( \cc O_{\PP(\E)}(-1) \otimes p^*\F^\vee\bigr) \times_B \prod_{i=j+1}^{r} \tot \bigl( \cc O_{\PP(\F^\vee)}(-1) \otimes p^*\E\bigr),
\]
where the products above are fiber products over $B$. Then $X_{T, {+}} = X_{T, r}$ and $X_{T, {-}} = X_{T, 0}$.

The equivariant cohomology ring $H^*_\torus(X_{T, j})$ is given by the quotient of
\[
H^*_\torus(B)[\z_i, \x_i]_{i=1}^n[\y_1, \dots, \y_r]
\]
by the relations
\begin{equation}\label{e:rels}
\left\{\prod_{i=1}^n (- \x_i - \y_k)\ \middle|\ 1 \leq k \leq j\right\} \bigcup \left\{\prod_{i=1}^n (\z_i + \y_k)\ \middle|\ j+1 \leq k \leq r\right\}.
\end{equation}

The $\torus$-fixed loci of $X_{T, \pm}$ are indexed by (not necessarily injective) functions $f\colon \{1, \dots, r\} \to \{1, \dots, n\}$. For such an $f$, the associated fixed locus $B_f$ is the subspace of $X_{T, +}$ such that, on the $i$th factor of
$\tot \bigl( \cc O_{\PP(\E)}(-1) \otimes p^*\F^\vee\bigr)$ in the fiber product, the homogeneous coordinates are all zero except for in the $f(i)$th factor of $\PP(\E)$.
In other words,
\[
B \cong B_f = \PP\bigl(\LB_{f(1)}\bigr) \times_B \cdots \times_B \PP\bigl(\LB_{f(r)}\bigr) \subset \PP(\E) \times_B \cdots \times_B \PP(\E) \subset X_{T, +}.
\]
An analogous description holds for $X_{T, -}$ and indeed all $X_{T, j}$.

We will make use of the following weak form of the abelian/non-abelian correspondence in cohomology.

\begin{Proposition}\label{p:pull}
The pullback
\[
H^*_\torus(X_\pm) \to \bigl(H^*_\torus\bigl( V^{s, \pm}(G)/T\bigr)\bigr)^W
\]
is an isomorphism, where $V^{s, \pm}(G)$ denotes the stable locus of $V$ with respect to $\det^{\pm}\colon G \to \CC^*$ and $W$ denotes the Weyl group. Furthermore, the pullback map
\[
H^*_\torus(X_{T, \pm}) \to H^*_\torus\bigl( V^{s, \pm}(G)/T\bigr)
\]
is surjective.
\end{Proposition}

See \cite[Proposition 1]{Brion} for details.
The proposition allows us to \emph{lift} cohomology classes in $X_{\pm}$ to $W$-invariant classes in $H^*_\torus(X_{T, \pm})$.

For future use we will use $F_\pm$ to denote the set of functions $f\colon \{1, \dots, r\} \to \{1, \dots, n\}$, i.e., $F_\pm$ indexes the set of fixed loci of $X_{T, \pm}$. We will use ${\rm In}_\pm$ to denote the set of injective functions $\{1, \dots, r\} \to \{1, \dots, n\}$, noting that ${\rm In}_\pm$ indexes the fixed loci of $V^{s, \pm}(G)/T.$

\subsection{Fourier--Mukai transform}\label{s:FM}
We can define a resolution of the birational map
$f\colon X_{-} \dashrightarrow X_{+}$
as follows. Let
\[
\tilde V = \mathcal{H}{\rm om}(\F, \cc O_B \otimes \CC^r) \times \mathcal{H}{\rm om}(\cc O_B \otimes \CC^r, \cc O_B \otimes \CC^r) \times \mathcal{H}{\rm om}(\cc O_B \otimes \CC^r, \E)
\]
and let $\tilde G = {\rm GL}_r \times {\rm GL}_r$. Consider the right action of $\tilde G$ on $\tilde V$ given by
\[
\bigl(\tilde M, \tilde D, \tilde N\bigr) \cdot (g_1, g_2)= \bigl( g_1^{-1}\tilde M , g_2^{-1}\tilde D g_1, \tilde Ng_2\bigr).
\]
Let $\tilde \theta \in \Hom\bigl(\tilde G, \CC^*\bigr)$ denote the character
\[
\tilde \theta(g_1, g_2) = \det(g_1)^{-1} \det(g_2).
\]
Define $\tilde X$ to be the GIT quotient,
$\tilde X = \tilde V \sslash_{\tilde \theta} \tilde G$.
The morphism $f_{-}\colon \tilde X \to X_{-}$ is induced by the morphism of stacks
$\bigl[\tilde V/\tilde G\bigr] \to [V/G]$
defined by the map
\begin{align*}
 \tilde V \to V, \qquad
 \bigl(\tilde M, \tilde D, \tilde N\bigr) \mapsto \bigl(\tilde M, \tilde N \tilde D\bigr)
\end{align*}
and the homomorphism
\begin{align*}
 \tilde G \to G ,\qquad
 (g_1, g_2) \mapsto g_1.
\end{align*}
There is a similar description of the map $f_{+}\colon \tilde X \to X_{+}.$
The $\torus$-action on $E$ and $F$ induces an action on $\tilde X$ for which the maps $f_\pm\colon \tilde X \to X_\pm$ are equivariant. We define
\begin{equation}\label{Kmap}
\FM\colon \ K^0_\torus(X_{-}) \to K^0_\torus(X_{+})
\end{equation}
to be $\FM = (f_{+})_* \circ f_{-}^*.$

\begin{Lemma}\label{l:FP} Let $H = \tot\bigl(F^\vee \otimes E\bigr) = \tot( \mathcal{H}{\rm om}(F, E) )$.
 The space $\tilde X$ is equal to the fiber product
 $X_+ \times_H X_-.$
\end{Lemma}
\begin{proof}
 It suffices to show that $\tilde X$ satisfies the universal property of fiber products. Suppose that $q\colon S \to B$ is a $B$-scheme. A morphism $S \to X_+$ over $B$ corresponds to a rank $r$ vector bundle $U_+ \to S$ together with linear maps $M_+\colon F \to U_+$ and $N_+\colon U_+ \to E$ such that $N_+$ is full rank. A morphism $S \to X_-$ over $B$ corresponds to a rank $r$ vector bundle $U_- \to S$ together with linear maps $M_-\colon F \to U_-$, $N_-\colon U_- \to E$ such that $M_-$ is full rank.

 On the other hand, a morphism $S \to \tilde X$ corresponds to a pair of rank $r$ vector bundles $U_{+/-} \to S$ (obtained via the induced map $S \to B{\rm GL}_r \times B{\rm GL}_r$) together with linear maps
 \smash{$\tilde M\colon F \to U_-$}, \smash{$\tilde D\colon U_- \to U_+$}, and \smash{$\tilde N\colon U_+ \to E,$} such that the first and last are both full rank.

 Suppose we have morphisms $S \to X_{+/-}$ which agree when composed with the natural maps $X_{+/-} \to H$.
 This yields two triples $(U_+, M_+, N_+)$ and $(U_-, M_-, N_-)$ such that $N_+ \circ M_+ = N_- \circ M_-$.
 We define a morphism \smash{$S \to \tilde X$} as follows.
 Note first that since $N_+$ is full rank, the condition $N_+ \circ M_+ = N_- \circ M_-$ implies that $N_-$ factors through a map \smash{$\tilde D\colon U_- \to U_+$}. Then define \smash{$\tilde M = M_-$} and \smash{$\tilde N = N_+.$} The triple $\tilde M$, $\tilde D$, $\tilde N$ defines a map $S \to \tilde X$ as in the previous paragraph.\looseness=-1

It is easy to check that the compositions $S \to \tilde X \to X_{+/-}$ agree with the original morphisms $S \to X_{+/-}$ from the previous paragraph. Thus \smash{$\tilde X$} satisfies the universal property as desired.
\end{proof}

In \cite[Theorem~D]{BLV}, it is proven that the
associated functor $\mathbb{FM}\colon D^b(X_-) \to D^b\bigl(X_+\bigr)$, between derived categories is an equivalence in the case that $B$ is a point.\footnote{In \cite{BLV}, the Fourier--Mukai kernel used is $\cc O_{X_+ \times_H X_-}$. This agrees with our kernel $\cc O_{\tilde X}$ by Lemma~\ref{l:FP}.} In~\cite{BCFMV}, it is shown that this equivalence may also be described using windows. The arguments naturally extend to our setting. We will use this fact only to conclude that the map~\eqref{Kmap} is an isomorphism and preserves the pairing.

We compute $\FM $ using localization as in \cite{Thomason}.
The fixed loci of $\tilde X$ are indexed by a pair of size-$r$ subsets of $\{1, \dots, n\}$,
corresponding to choosing $r \times r$ minors of both $\mathcal{H}{\rm om}(\cc O_B \otimes \CC^r, \E)$ and $\mathcal{H}{\rm om}(\F, \cc O_B \otimes \CC^r)$, i.e., they are indexed by elements of $D_- \times D_+$.

Let $R_1$, $R_2$ be the representations of $\tilde G$ obtained by composing the projection $\pi_1, \pi_2\colon \tilde G \to G$ to $G$ with the right representation $R$ of $G$. Again, we abuse notation and also use \smash{$R_{1/2}$} to denote the induced vector bundles on $\tilde X$.
With these conventions, we observe that \smash{$f_{-/+}^*(R) = R_{1/2}$}.

A simple generalization of Lemma~\ref{l:rest}, gives the following.
\begin{Lemma}\label{l:rest2}
For $(\delta^-, \delta^+) \in D_{-} \times D_{+}$,
let $B_{(\delta^-, \delta^+)} $ denote the associated fixed locus, isomorphic to $B$. Then $f_{+}$ $($resp.\ $f_-)$ maps $B_{(\delta^-, \delta^+)} $ isomorphically onto $B_{\delta^+}$ $($resp.\ $B_{\delta^-})$. Furthermore,
$(R_1)|_{B_{(\delta^-, \delta^+)}} \cong \bigoplus_{i=1}^r \MB_{\delta^-_i}^\vee$ and $(R_2)|_{B_{(\delta^-, \delta^+)}} \cong \bigoplus_{j=1}^r \LB_{\delta^+_j}^\vee$.
\end{Lemma}

In order to ease notation, in this section we will also use $L_i$ or $M_j$ to denote the pullbacks of these respective line bundles to $X_+$, $X_-$, and $\tilde X.$
With this convention, for $\delta^- \in D_-$, let $e_{\delta^-}$ denote the class{\samepage
\begin{align}
e_{\delta^-} = \prod_{j_- \notin \delta^-} \wedge^\bullet \mathcal{H}{\rm om}(\MB_{j_-}, \cc O_B \otimes \CC^r)^\vee
\cong \prod_{j_- \notin \delta^-} \wedge^\bullet \bigl(\MB_{j_-}^\vee \otimes R^\vee\bigr)^\vee \label{edm}
\end{align}
in $K^0_\torus(X_-)$, where for $A \in K^0_\torus(X_-)$, we denote $\wedge^\bullet A := \bigoplus_{i=0}^{\op{rk}(A)} (-1)^i \wedge^i A$.}

\begin{Corollary}\label{c:gen}
The $\torus$-localized equivariant $K$-theory of $X_-$ is spanned by the classes $\bigl\{\pi_-^*(A) \otimes e_{\delta^-} | A \in K^0(B), \delta^- \in D_{-}\bigr\}$.
\end{Corollary}

\begin{proof}
By Lemma~\ref{l:rest}, one checks that the class $\pi_-^*(A) \otimes e_{\delta^-}$ restricts to zero on \smash{$B_{ \delta_0^-}$} for $\delta_0^- \neq \delta^-$. Furthermore, generalizing the argument of \cite[Section~5.1]{CKM14} to the relative setting, we see that the relative tangent bundle of $\Gr\bigl(r, \F^\vee\bigr)$ is obtained from the following \emph{relative Euler sequence} for $\Gr\bigl(r, F^\vee\bigr)$:
\begin{equation}\label{e:Euler} 0 \rightarrow S_{-} \otimes S_{-}^\vee \rightarrow S_{-}^\vee \otimes F^\vee \rightarrow T\Gr\bigl(r, \F^\vee\bigr)/ B \rightarrow 0.\end{equation} Therefore, $\pi_-^*(A) \otimes e_{\delta^-}$ restricts to \smash{$A \otimes \wedge^\bullet N_{B_{\delta^-}|\Gr(r, \F^\vee)}^\vee$} where \smash{$\wedge^\bullet N_{B_{\delta^-}|\Gr(r, \F^\vee)}^\vee$} is an invertible class on $B_{\delta^-}$. The result then follows from the localization theorem in $K$-theory as in \mbox{\cite[Section~2.3]{Oko15}}.
\end{proof}

\begin{Proposition}
We have
\begin{equation}\label{fmedm}
\FM(e_{\delta^-}) = \prod_{j_- \notin \delta^-} \wedge^\bullet \bigl(\MB_{j_-}^\vee \otimes R^\vee\bigr)^\vee.
\end{equation}
\end{Proposition}

\begin{Remark}
 Note that although the right-hand side of the equation above is the same formal expression as the right-hand side of \eqref{edm}, the right-hand side of \eqref{fmedm} is to be interpreted as a~class in $K^0_\torus\bigl(X_+\bigr)$.
\end{Remark}

\begin{proof}
Note first that,
in $K^0_\torus\bigl(\tilde X\bigr)$, the class
\[
\prod_{ i_{-} \notin \delta^-} \wedge^\bullet \bigl(\MB_{i_{-}}^\vee \otimes R_1^\vee\bigr)^\vee \cdot \prod_{j_{+} \notin \delta^+} \wedge^\bullet ( \LB_{j_+} \otimes R_2)^\vee \cdot \wedge^\bullet \bigl(R_2^\vee \otimes R_1\bigr)^\vee
\]
restricts to \smash{$\wedge^\bullet N_{B_{\delta^-, \delta^+}|\tilde X}^\vee$} on
\smash{$B_{(\delta^-, \delta^+)}$}, and restricts to zero on \smash{$B_{(\delta^-_1, \delta^+_1 )}$} for \smash{$\bigl(\delta^-_1, \delta^+_1\bigr) \neq (\delta^-, \delta^+)$}.
It follows that
\begin{gather*}
f_{-}^*(e_{\delta^-}) = \prod_{j_- \notin \delta^-} \wedge^\bullet \bigl(\MB_{j_-}^\vee \otimes R_1^\vee\bigr)^\vee \\
\hphantom{f_{-}^*(e_{\delta^-})}{}
= \sum_{\delta^+ \in D_{+}} (i_{\delta^-, \delta^+})_*i_{\delta^-, \delta^+}^* \\
\hphantom{f_{-}^*(e_{\delta^-}) =\sum_{\delta^+ \in D_{+}}}{}
 \times\Biggl(\frac{\prod_{j_- \notin \delta^-} \wedge^\bullet \bigl(\MB_{j_-}^\vee \otimes R_1^\vee\bigr)^\vee}
{ \prod_{ i_{-} \notin \delta^-} \wedge^\bullet \bigl(\MB_{i_{-}}^\vee \otimes R_1^\vee\bigr)^\vee \cdot \prod_{j_{+} \notin \delta^+} \wedge^\bullet ( \LB_{j_+} \otimes R_2)^\vee \cdot \wedge^\bullet \bigl(R_2^\vee \otimes R_1\bigr)^\vee}\Biggr)
\\
\hphantom{f_{-}^*(e_{\delta^-})}{}
= \sum_{\delta^+ \in D_{+}} (i_{\delta^-, \delta^+})_*i_{\delta^-, \delta^+}^* \Biggl(\frac{1}
{ \prod_{j_{+} \notin \delta^+} \wedge^\bullet ( \LB_{j_+} \otimes R_2)^\vee \cdot \wedge^\bullet \bigl(R_2^\vee \otimes R_1\bigr)^\vee}\Biggr)
\\
\hphantom{f_{-}^*(e_{\delta^-})}{}
= \sum_{\delta^+ \in D_{+}} (i_{\delta^-, \delta^+})_*i_{\delta^-, \delta^+}^* \Biggl(\frac{\prod_{j_- \notin \delta^-} \wedge^\bullet \bigl(\MB_{j_-}^\vee \otimes R_2^\vee\bigr)^\vee }
{ \prod_{j_{+} \notin \delta^+} \wedge^\bullet ( \LB_{j_+} \otimes R_2)^\vee \cdot \prod_{1 \leq j \leq n} \wedge^\bullet \bigl(\MB_{j}^\vee \otimes R_2^\vee\bigr)^\vee }\Biggr),
\end{gather*}
where the second equality is by the localization theorem applied to $K^0_\torus\bigl(\tilde X\bigr)$, the fourth is Lemma~\ref{l:rest2}.
Pushing forward via $f_{+}$, we obtain

\begin{align*}
(f_{+})_* f_{-}^*(e_{\delta^-})
={}& \sum_{\delta^+ \in D_{+}} (i_{\delta^+})_*i_{\delta^+}^* \Biggl(\frac{\prod_{j_- \notin \delta^-} \wedge^\bullet \bigl(\MB_{j_-}^\vee \otimes R^\vee\bigr)^\vee }
{ \prod_{j_{+} \notin \delta^+} \wedge^\bullet ( \LB_{j_+} \otimes R)^\vee \cdot \prod_{1 \leq j \leq n} \wedge^\bullet \bigl(\MB_{j}^\vee \otimes R^\vee\bigr)^\vee }\Biggr) \\
={}& \prod_{j_- \notin \delta^-} \wedge^\bullet \bigl(\MB_{j_-}^\vee \otimes R^\vee\bigr)^\vee,
\end{align*}
where the last equality is the localization theorem applied to $K^0_\torus(X_{+})$.
\end{proof}

For $(\delta^-, \delta^+) \in D_{-} \times D_{+}$, define
\begin{align}\label{CD}
C_{\delta^-, \delta^+}
={}& \prod_{i=1}^r {\rm e}^{(n-r)(x_{\delta^+_i} - z_{\delta^-_i})/2} \prod_{j_- \notin \delta^-} \frac{\sin\bigl(\bigl(\x_{\delta^+_i} - \z_{j_-}\bigr)/2 {\rm i}\bigr)}{\sin\bigl(\bigl(\z_{\delta^-_i} - \z_{j_-}\bigr)/2 {\rm i}\bigr)},
\end{align}
where ${\rm i} = \sqrt{-1}$.

\begin{Definition}
\label{d:UH0}
Define $\UU_H\colon H^*_\torus(X_{-})\otimes_{R_\torus} S_\torus \to H^*_\torus(X_{+})\otimes_{R_\torus} S_\torus$ to be the linear map defined by sending
$\alpha/e_\torus(N_{\delta^-}) \in H^*_\torus(B_{\delta^-})\otimes_{R_\torus} S_\torus$ to
\begin{equation}\label{e:UH}
\sum_{\delta^+ \in D_+} C_{\delta^-, \delta^+}
\frac{ \phi_{\delta^+, \delta^-}(\alpha)}{e_\torus(N_{\delta^+})},
\end{equation}
where $\phi_{\delta^+, \delta^-}\colon H^*_\torus(B_{\delta^-}) \to H^*_\torus(B_{\delta^+})$ is the canonical isomorphism induced by the projections $\pi_{\delta^+}\colon B_{\delta^+} \to B$ and $\pi_{\delta^-}\colon B_{\delta^-} \to B$.
\end{Definition}

\begin{Proposition}\label{p:Kcom}
The following diagram commutes:
\[
\begin{tikzcd}
K^0_\torus(X_{-}) \ar[r, "\FM"] \ar[d, "\ch^\torus"] & K^0_\torus(X_{+}) \ar[d, "\ch^\torus"] \\ H^*_\torus(X_{-})\otimes_{R_\torus} \hat S_\torus \ar[r, "\UU_H"] & H^*_\torus(X_{+})\otimes_{R_\torus} \hat S_\torus.
\end{tikzcd}
\]
\end{Proposition}

\begin{proof}
We will prove
$\ch^\torus(\FM(e_{\delta^-})) = \UU_H \bigl( \ch^\torus(e_{\delta^-})\bigr)$.
 The general statement then follows by Corollary~\ref{c:gen} and the projection formula.
Observe that
\begin{align*}
\ch^\torus( i_{\delta^-}^* e_{\delta^-}) &=\ch^\torus\biggl( i_{\delta^-}^* \prod_{j_- \notin \delta^-} \wedge^\bullet \bigl(\MB_{j_-}^\vee \otimes R^\vee\bigr)^\vee \biggr)
\\
& =\ch^\torus\Biggl( \prod_{i=1}^r \prod_{j_- \notin \delta^-} \wedge^\bullet \bigl(\MB_{j_-}^\vee \otimes \MB_{\delta^-_i} \bigr)^\vee \Biggr)
=
\prod_{i=1}^r \prod_{j_- \notin \delta^-} \Bigl( 1 - {\rm e}^{\z_{\delta^-_i} - \z_{j_-}}\Bigr)
\end{align*}
and
\begin{align*}
\ch^\torus( i_{\delta^+}^* \FM(e_{\delta^-} ))& = \ch^\torus i_{\delta^+}^*\biggl( \prod_{j_- \notin \delta^-} \wedge^\bullet \bigl(\MB_{j_-}^\vee \otimes R^\vee\bigr)^\vee
\biggr) \\
& =
\ch^\torus\Biggl( \prod_{i=1}^r \prod_{j_- \notin \delta^-} \wedge^\bullet \bigl(\MB_{j_-}^\vee \otimes \LB_{\delta^+_i} \bigr)^\vee \Biggr)
=
\prod_{i=1}^r \prod_{j_- \notin \delta^-} \Bigl( 1 - {\rm e}^{\x_{\delta^+_i} - \z_{j_-}}\Bigr).
\end{align*}
The ratio of the two expressions above, after identifying $B_{\delta^-}$ and $B_{\delta^+}$, is
\begin{align*}
& \prod_{i=1}^r \prod_{j_- \notin \delta^-} \frac{\Bigl( 1 - {\rm e}^{\x_{\delta^+_i} - \z_{j_-}}\Bigr)}{\Bigl( 1 - {\rm e}^{\z_{\delta^-_i} - \z_{j_-}}\Bigr)} \\
&\qquad{}= \prod_{i=1}^r \prod_{j_- \notin \delta^-} \frac{ {\rm e}^{(\x_{\delta^+_i} - \z_{j_-})/2}\Bigl( {\rm e}^{-(\x_{\delta^+_i} - \z_{j_-})/2} - {\rm e}^{(\x_{\delta^+_i} - \z_{j_-})/2}\Bigr)}{ {\rm e}^{(\z_{\delta^-_i} - \z_{j_-})/2}\Bigl( {\rm e}^{-(\z_{\delta^-_i} - \z_{j_-})/2} - {\rm e}^{(\z_{\delta^-_i} - \z_{j_-})/2}\Bigr)}
= C_{\delta^-, \delta^+}.
\end{align*}
By the localization theorem in cohomology (see \cite{AB}), it follows that the diagram commutes.
\end{proof}

\begin{Corollary}\label{cv}
The map $\UU_H$ is invertible and has a well-defined non-equivariant limit. Furthermore, we have
\[
{\pi_-}_*\bigl( \alpha \cup \op{Td}_{X_-/B} \bigr) = {\pi_+}_*\bigl( \UU_H(\alpha) \cup \op{Td}_{X_+/B} \bigr).
\]
\end{Corollary}
\begin{proof}

Let $\ol X_+ $ and $ \ol X_-$ denote the GIT quotients from Section~\ref{gs} in the special case that the base is a point, so $\ol X_+ \cong \ol X_- = \tot(S)$, where $S$ is the tautological bundle over $\Gr(r, \CC^n)$.
Let
\[
\ol \FM\colon \ \hat K^0_\torus\bigl(\ol X_-\bigr) \to \hat K^0_\torus\bigl(\ol X_+\bigr)
\]
denote the Fourier--Mukai transform in this case, where \smash{$\hat K^0_\torus\bigl(\ol X_\pm\bigr)$} denotes the completion of \smash{$K^0_\torus\bigl(\ol X_\pm\bigr)$} along the augmentation ideal (see~\cite{EG}).

By definition of the $\torus$ action on $B$, $E$ and $F$, the (completion of the) equivariant cohomology of $B$ is given by
\[
\hat H^*_\torus(B) = H^*(B)[\![\z_i, \x_i]\!]_{i=1}^n,
\]
where recall that $\x_i = c_1^\torus\bigl(\LB_i^\vee\bigr)$ and $\z_i = c_1^\torus\bigl(\MB_i^\vee\bigr)$.
The equivariant line bundles $\LB_i \to B$ and $\MB_i \to B$ for $1 \leq i \leq n$ determine a map
\[v\colon \ \big[B/(\CC^*)^{2n}\big] \to \big[\{pt\}/(\CC^*)^{2n}\big].
\]
This induces the following commutative diagram of stacks:
\[
 \begin{tikzcd}
 \bigl[X_{\pm}/(\CC^*)^{2n}\bigr] \ar[d, "\pi_\pm"] \ar[r, "\tilde v_\pm"] & \bigl[\ol X_{\pm}/(\CC^*)^{2n}\bigr] \ar[d, "\ol \pi_\pm"] \\
 \bigl[B/(\CC^*)^{2n}\bigr] \ar[r, "v"] & \bigl[\{pt\}/(\CC^*)^{2n}\bigr].
 \end{tikzcd}
\]
Taking cohomology, we obtain
\begin{equation}\label{gdiag}
 \begin{tikzcd}
 H^*(B)[\![\z_i, \x_i]\!]_{i=1}^n[\sigma_1, \dots, \sigma_r]/\hat I_\pm \ar[d, "{\pi_\pm}_*"] & \ar[l, "\tilde v_\pm^*"] \CC[\![\z_i, \x_i]\!]_{i=1}^n[\sigma_1, \dots, \sigma_r]/\hat I_\pm \ar[d, "{\ol \pi_\pm}_*"] \\
 H^*(B)[\![\z_i, \x_i]\!]_{i=1}^n & \ar[l, "v^*"] \CC[\![\z_i, \x_i]\!]_{i=1}^n.
\end{tikzcd}
\end{equation}
This diagram commutes. In fact, after identifying
\[
\hat H^*_\torus(X_{\pm}) := \prod_{i=0}^\infty H^i_\torus(X_{\pm})
\]
with $H^*(B) \otimes_\CC \CC[\![\z_i, \x_i]\!]_{i=1}^n[\sigma_1, \dots, \sigma_r]/\hat I_\pm,$ we have the identification
\[
{\pi_\pm}_* = \op{id}_B \otimes_\CC {\ol \pi_{\pm}}_*,
\]
where $\op{id}_B$ denotes the identity map on $H^*_\torus(B)$.

After taking completions with respect to the augmentation ideal, the Chern map
\[
\hat{\ch}^\torus\colon \ \hat K^0_\torus\bigl(\ol X_\pm\bigr) \to \hat H^*_\torus\bigl(\ol X_\pm\bigr)
\]
is an isomorphism \big(see \cite{EG}; here we use $\hat{\ch}^\torus$ to denote the induced map defined on the completion \smash{$\hat K^0_\torus\bigl(\ol X_\pm\bigr)$}\big). Because it is induced by a derived equivalence, the map $\ol \FM$ is also an isomorphism.
Proposition~\ref{p:Kcom} implies that \smash{$\ol \UU_H := \hat{\ch}^\torus \circ \ol \FM \circ \bigl(\hat{\ch}^\torus\bigr)^{-1}$}
is also an
isomorphism, and may be defined before localizing with respect to $\torus$. In other words, $\ol \UU_H$ gives a well-defined map
\[
\ol \UU_H\colon \ \CC[\![\z_i, \x_i]\!]_{i=1}^n[\sigma_1, \dots, \sigma_r]/\hat I_- \to \CC[\![\z_i, \x_i]\!]_{i=1}^n[\sigma_1, \dots, \sigma_r]/\hat I_+.
\]

Again using the identification
$ \hat H^*_\torus(X_{\pm}) = H^*(B) \otimes_\CC \CC[\![\z_i, \x_i]\!]_{i=1}^n[\sigma_1, \dots, \sigma_r]/\hat I_\pm$
together with the explicit formula for $\UU_H$ in \eqref{e:UH}, we also see that $\UU_H = \op{id}_B \otimes_\CC \ol \UU_H$ which implies the first statement.

Given $\ol \alpha \in \hat H^*_\torus\bigl(\ol X_-\bigr)$, let $A = \bigl(\ch^\torus\bigr)^{-1}(\ol \alpha) \in \hat K^0_\torus\bigl(\ol X_-\bigr)$. Applying the equivariant Hirzebruch--Riemann--Roch theorem (see \eqref{e:HRR} below) twice, we see that for \smash{$A \in \hat K^0_\torus(X_-)$},
\begin{align}
{\ol \pi_-}_* \bigl(\ol \alpha \op{Td}_{\ol X_-} \bigr)& = \chi(A)
 = \chi\bigl(\ol \FM (A)\bigr)
 = {\ol \pi_+}_* \bigl(\ch^\torus\bigl(\ol \FM (A)\bigr) \op{Td}_{\ol X_+}\bigr) \nonumber \\
 & = {\ol \pi_+}_* \bigl(\UU_H (\ol \alpha) \op{Td}_{\ol X_+}\bigr). \label{UC}
 \end{align}

The pullback $\tilde v_\pm^*\bigl(\op{Td}_{\ol X_\pm}\bigr)$ is the relative Todd class $\op{Td}_{X_\pm}$. This fact together with \eqref{UC} allows us to upgrade \eqref{gdiag} to the following commutative diagram:
\begin{equation}\label{vgdiag}
 \begin{tikzcd}
 \CC[\![\z_i, \x_i]\!]_{i=1}^n[\sigma_k]/\hat I_- \ar[rr, "\ol \UU_H"] \ar[dr, swap, "{\ol \pi_-}_*( \op{Td}_{\ol X_-} \cup - )"] \ar[dd, swap, "\tilde v_-^*"] & & \CC[\![\z_i, \x_i]\!]_{i=1}^n[\sigma_k]/\hat I_+ \ar[dl, "{\ol \pi_+}_*( \op{Td}_{\ol X_+} \cup - )"] \ar[dd, "\tilde v_+^*"]\\
 & \CC[\![\z_i, \x_i]\!]_{i=1}^n \ar[dd, near start, crossing over, "v^*"]
 & \\
 H^*(B)[\![\z_i, \x_i]\!]_{i=1}^n[\sigma_k]/\hat I_- \ar[rr, near start, "\UU_H"] \ar[dr, swap, "{\pi_-}_*(\op{Td}_{X_-/B} \cup - )"] && H^*(B)[\![\z_i, \x_i]\!]_{i=1}^n[\sigma_k]/\hat I_+ \ar[dl, "{\pi_+}_*( \op{Td}_{X_+/B} \cup - )"] \\
 &H^*(B)[\![\z_i, \x_i]\!]_{i=1}^n \ar[from=uu, near start, crossing over, "v^*"]&
 \end{tikzcd}
\end{equation}
The second claim follows.
\end{proof}

\section[Gromov-Witten theory and I-functions]{Gromov--Witten theory and $\boldsymbol{I}$-functions}\label{s:I}

Let $X$ be a smooth semi-projective variety. In what follows, $X$ will also be equipped with an action by a torus $\torus$. Following \cite{CIJS}, we assume that all $\torus$-weights appearing in the representation $H^0(X, \cc O_X)$ lie in a strictly convex cone, and that $H^0(X, \cc O_X)^\torus = \CC.$ This implies in particular that $X^\torus$ is compact.

Following standard practice, we define integrals over $X$ via localization. For $\gamma \in H^*_\torus(X)$, define
\[
\int_X \gamma := \sum_F \int_F \frac{(\gamma)|_F}{e_\torus\bigl(N_{F|X}\bigr)} \in S_\torus,
\]
where the sum is over all connected components $F$ in $X^\torus$ and the class $e_\torus( - )$ denotes the $\torus$-equivariant Euler class. Define a pairing $\br{-,-}_X$ on $H^*_\torus(X) \otimes_{R_\torus} S_\torus$ by
\[
\br{\alpha, \beta}_X := \sum_F \int_F \frac{(\alpha\cup \beta)|_F}{e_\torus\bigl(N_{F|X}\bigr)}.
\]

For a ring $R$, let
\[
R[\![Q]\!] = \biggl\{ \sum_{d\in \eff(X)} a_d Q^d \mid a_d \in R \biggr\}
\]
denote the \emph{Novikov ring} of $X$.
Let
\[
\cc H_X = H^*_\torus(X) \otimes_{R_\torus} S_\torus\bigl(\!\bigl(z^{-1}\bigr)\!\bigr)[\![Q]\!]
\]
denote Givental's symplectic vector space, equipped with the symplectic form
\begin{equation*}
\Omega(f, g) = \op{Res}_{z = 0}\br{f(-z), g(z)}_X {\rm d}z.
\end{equation*}

This defines a polarization
\[
\cc H_X = \cc H_X^+ \oplus \cc H_X^-
\]
with $\cc H_X^+ = H^*_\torus(X) \otimes_{R_\torus} S_\torus[z][\![Q]\!]$ and
$\cc H_X^- = z^{-1} H^*_\torus(X) \otimes_{R_\torus} S_\torus \bigl[\!\bigl[z^{-1}\bigr]\!\bigr] [\![Q]\!]$.

For $\alpha_1, \dots, \alpha_m \in H^*_\torus(X)$, $k_1, \dots, k_m \geq 0$, and $d \in \eff(X)$, define the degree-$d$ genus-0 Gromov--Witten invariant of $X$:
\[
\big\langle\alpha_1 \psi_1^{k_1}, \dots, \alpha_m \psi_m^{k_m}\bigr\rangle_{0, m ,d}^X = \int_{[\cc M_{0, m}(X, d)]^{vir}} \prod_{i=1}^m ev_i^*(\alpha_i) \psi_i^{k_i}.
\]
If $X$ is not compact, the integral on the right-hand side is defined after localization.
\begin{Definition}
Define the $\torus$-equivariant \emph{genus-zero descendent potential} to be
\begin{equation*}
\cc F_X(\bt(z)) = \sum_{d \in \eff(X)} \sum_{m \geq 0} \frac{Q^d}{m!} \br{\bt(\psi_1), \dots, \bt(\psi_m)}_{0, m, d}^X,
\end{equation*}
for $\bt(z) = \sum_{a=0}^s t_a z^a \in H^*_\torus(X)[z]$. Extending linearly, $\cc F_X$ defines a function on $\cc H_X^+$.
\end{Definition}

\begin{Definition}
After applying the \emph{dilaton shift} $\sum_{a=0}^s t_a z^a = z1 + \sum_{a=0}^s q_a z^a$, we define Givental's \emph{overruled Lagrangian cone} $\cc L_X$ to be the graph of $\cc F_X$ with respect to the above polarization:
\[
\cc L_X = \bigl\{ (p, q) \in T^*\cc H_X^+ \mid p = d_q \cc F_X\bigr\}.
\]
\end{Definition}
Let $\{\phi_i\}_{i \in I}$ be a basis for $H^*_\torus(X)$ and let $\{\phi^i\}_{i \in I}$ denote the dual basis.
Points of $\cc L_X$ are of the form
\[
-z 1 + \bt(z) + \sum_{d \in \eff(X)} \sum_{m \geq 0} \frac{Q^d}{m!} \br{\frac{\phi_i}{-z - \psi_1}, \bt(\psi_2), \dots, \bt(\psi_{m+1})}_{0, m+1, d}^X \phi^i.
\]
\begin{Definition}
The $J$-function of $X$ is the slice of $\cc L_X$ obtained by intersecting with $-z 1 + t + \cc H_X^-$ for $t = \sum_{i \in I} t^i \phi_i \in H^*_\torus(X)$:
\[
J^X(t, z) = -z 1 + t+ \sum_{d \in \eff(X)} \sum_{m \geq 0} \frac{Q^d}{m!} \sum_{i \in I}\br{\frac{\phi_i}{-z - \psi_1}, t, \dots, t}_{0, m+1, d}^X \phi^i.
\]

\end{Definition}
The function $J^X(t, z)$ determines the entire Lagrangian cone. Via the string equation, dilaton equation, and topological recursion relations, any point on $\cc L_X$ may be expressed as
\[
J^X(\tau , z) + \sum_{i \in I} C^i(z) \frac{\partial}{\partial{t^i}} J^X(t, z)|_{t=\tau},
\]
where $C^i(z) \in S_\torus[z][\![Q]\!]$ are \emph{polynomials} in $z$.

\subsection{Integral lattice}
As predicted by mirror symmetry, the genus zero Gromov--Witten theory ($A$-model) of a space~$X$ should be equivalent to the Gauss--Manin connection ($B$-model) of a family of mirror manifolds~$\check{X_t}$. The integral cohomology~$H^*\bigl(\check{X_t}, \ZZ\bigr)$ naturally endows the Gauss--Manin connection with an integral lattice of flat sections. Thus by mirror symmetry, one should expect an integral structure to appear in Gromov--Witten theory as well. Such a structure was defined by Iritani in \cite{IriInt} where it was verified to be compatible with mirror symmetry in the toric setting.

In this section, we describe Iritani's integral structure on $\tilde H_X$. We define a lattice in $\tilde H_X$ as the image of $K_\torus(X)$ under a certain modification of the Chern character. A key ingredient is a characteristic class known as the $\hat\Gamma$ class, which plays the role of a square root of the Todd class. This integral structure will allow us to compare the symplectic transformation $\mathbb{U}$ with the Fourier--Mukai transform $\FM$ from Section~\ref{s:FM}.

For $\phi \in H^{2p}_\torus(X)$ a homogeneous element, let
\[
\deg(\phi) = 2p
\]
denote the real degree.
We will use the following three operators heavily in what follows.

First, we define $\mu$ to be the operator
\[
\mu(\phi) = \biggl(\frac{\deg(\phi)}{2} - \frac{\dim(X)}{2}\biggr)\phi,
\]
where $\dim$ denotes the dimension over $\CC$.

Secondly, we define $\deg_0\colon H^*_\torus(X) \to H^*_\torus(X)$ be the degree operator, defined by
\[
\deg_0(\phi) = \deg(\phi)\cdot \phi.
\]

Finally, we define $\rho$ to be the operator of multiplication by $c_1(T_X)$.

For $A \to X$ a vector bundle, define the multiplicative class $\hat \Gamma(A)$ to be
\[
\hat \Gamma(A) = \prod_{i=1}^s \Gamma(1 + a_i),
\]
where $a_1, \dots, a_s$ are the Chern roots of $A$, and $\Gamma(z)$ is the Gamma function expanded as a power series at $z = 0$. Denote $\hat \Gamma(TX)$ by $\hat \Gamma_{X}$.

Define the map
\begin{align*}
\Psi\colon \ K^0_\torus(X) \to \widetilde{\cc H}_X, \qquad
A \mapsto z^{- \mu} z^{\rho} \bigl( \hat \Gamma_X \cup (2 \pi {\rm i})^{\op{deg}_0/2} \ch^\torus(A) \bigr),
\end{align*}
where
\[
\widetilde{\cc H}_X = \cc H_X[\log z]\big[z^{-1/2}\big].
\]

The image of $\Psi$ defines a lattice in $\widetilde{\cc H}_X$. For future use, we will also define the map
\begin{align}
\psi\colon \ \hat H^*_\torus(X) \to \widetilde{\cc H}_X, \qquad
\phi \mapsto z^{- \mu} z^{\rho} \bigl( \hat \Gamma_X \cup (2 \pi {\rm i})^{\op{deg}_0/2} \phi \bigr)\label{psi}
\end{align}
so that $\Psi = \psi \circ \ch^\torus$.

Define the equivariant Euler characteristic to be
\[
\chi(A) := \sum_{i=0}^{\dim(X)} (-1)^i \ch^\torus\bigl(H^i(X, A)\bigr),
\]
where, as before, we use the equivariant Chern character with respect to $\torus$.
If $X$ is compact, then by \cite{EG, Thomason} we have the following equivariant Hirzebruch--Riemann--Roch theorem
\begin{equation}\label{e:HRR}\chi(A) = \int_X \ch^\torus(A) \op{Td}_X,\end{equation}
as an equality in $\hat S_\torus$ (as defined in Section~\ref{s:coh}).
 When $X$ is not compact the equivariant Hirzebruch--Riemann--Roch theorem is not known in general, but the following simple case follows easily from the projection formula.
Suppose that $Z$ is a smooth projective variety equipped with an action of $\torus$, that $V \to Z$ is an equivariant vector bundle, and that $X$ is the total space of~$V$. Suppose further that the $\torus$-fixed locus of $X$ is contained in $Z$.
Then the projection formula in cohomology implies
\[
\int_X \alpha = \int_Z s^*(\alpha)/e_\torus(V),
\]
where $s\colon Z \to V$ denotes the zero section.
On the other hand, the projection formula in $K$-theory implies
\[
\chi(A) = \chi\bigl(s^*(A)/\wedge^\bullet\bigl(V^\vee\bigr)\bigr).
\]
Applying the Hirzebruch--Riemann--Roch theorem \eqref{e:HRR} to $Z$ then implies \eqref{e:HRR} for $X$.

Following \cite{CIJ}, consider a modified Euler characteristic $\chi_z$ given by
\[
\chi_z(C, D) := \biggl(\frac{2 \pi {\rm i}}{ z}\biggr)^{\lambda \partial_\lambda} \chi(C, D),
\]
where $\lambda \partial_\lambda := \sum_i \lambda_i \frac{\partial}{\partial \lambda_i}$ for some basis $\lambda_1, \dots, \lambda_s$ of $H^2(B\torus)$.

\begin{Lemma}[{\cite[Proposition~2.10]{IriInt}}]
The map $\Psi$ relates the pairings in $K$-theory and cohomology as follows:
\[
\chi_z(C, D) = \frac{1}{ (2 \pi)^{\dim X}}\br{\Psi(C)(-z), \Psi(D)(z)}_X.
\]
\end{Lemma}
\begin{proof}
Following \cite[Proposition~3.2]{CIJ}, we use the following equalities:
\[
 z^{-\lambda {\rm d}\lambda} \br{\alpha,\beta}_X=\br{z^{-\mu}\alpha,z^{-\mu}\beta}_X$, $\br{z^{\rho}\alpha,z^\rho \beta}_X=\br{\alpha,\beta}_X
 \qquad \text{and}\qquad {\rm e}^{2\pi{\rm i}\mu}\rho=-\rho {\rm e}^{2\pi{\rm i}\mu}.
 \]
The proof then follows from the equivariant Hirzebruch--Riemann--Roch theorem,
\begin{align*}
&\chi_z(C, D)
 = z^{- \lambda \partial \lambda } \frac{1}{(2 \pi {\rm i})^{\dim X}} \int_X (2 \pi {\rm i})^{\op{deg}_0/2} \ch^\torus\bigl(C^\vee \otimes D\bigr) \op{Td}_X \\
&\qquad{} = z^{-\lambda \partial \lambda } \frac{1}{(2 \pi {\rm i})^{\dim X}} \int_X
 {\rm e}^{ \pi {\rm i} \rho} \hat \Gamma_X^\vee \hat \Gamma_X
(2 \pi {\rm i})^{\op{deg}_0/2} \ch^\torus\bigl(C^\vee \otimes D\bigr) \\
&\qquad{} = \frac{z^{-\lambda \partial \lambda }}{ (2 \pi)^{\dim X}} \int_X {\rm e}^{-\pi {\rm i} \dim(X)/2}
 {\rm e}^{\pi {\rm i} \rho} \hat \Gamma_X^\vee \hat \Gamma_X (2 \pi {\rm i})^{\op{deg}_0/2} \ch^\torus\bigl(C^\vee \otimes D\bigr)
 \\
&\qquad{}=
\frac{z^{-\lambda \partial \lambda }}{ (2 \pi)^{\dim X}}\bigl\langle {\rm e}^{\pi {\rm i} \rho} {\rm e}^{\pi {\rm i} \mu} \hat \Gamma_X (2 \pi {\rm i})^{\op{deg}_0/2}\ch^\torus(C), \hat \Gamma_X (2 \pi {\rm i})^{\op{deg}_0/2} \ch^\torus(D)\bigr\rangle_X \\
&\qquad{}=
\frac{1 }{ (2 \pi)^{\dim X}}\bigl\langle z^{- \mu} z^{- \rho} {\rm e}^{\pi {\rm i} \rho} {\rm e}^{\pi {\rm i} \mu} \hat \Gamma_X (2 \pi {\rm i})^{\op{deg}_0/2}\ch^\torus(C), z^{- \mu} z^{\rho} \hat \Gamma_X (2 \pi {\rm i})^{\op{deg}_0/2} \ch^\torus(D)\bigr\rangle_X \\
&\qquad{}=
\frac{1}{ (2 \pi)^{\dim X}}\bigl\langle\bigl( {\rm e}^{-\pi {\rm i}} z\bigr)^{- \mu} \bigl( {\rm e}^{-\pi {\rm i}}z\bigr)^\rho \hat \Gamma_X (2 \pi {\rm i})^{\op{deg}_0/2}\ch^\torus(C), z^{- \mu} z^{\rho} \hat \Gamma_X (2 \pi {\rm i})^{\op{deg}_0/2} \ch^\torus(D)\bigr\rangle_X \\
&\qquad{}= \frac{1}{ (2 \pi)^{\dim X}}\bigl\langle\Psi(C)\bigl( {\rm e}^{-\pi {\rm i}} z\bigr), \Psi(D)(z)\bigr\rangle_X.
\tag*{\qed}
\end{align*}
\renewcommand{\qed}{}
\end{proof}

Let us return now to the specific geometry of the relative Grassmann flop of the previous section. Let $X_+$ and $X_-$ be as in Section~\ref{gs}. Let $\psi_{\pm}$ denote \eqref{psi} for these particular spaces.
\begin{Definition}\label{d:UH}
Define $\UU\colon \cc H_{X_-} \to \cc H_{X_+}$ to be
$\UU = \psi_{+} \circ \UU_H \circ \psi_{-}^{-1}$,
 where $\UU_H$ is the map in~\eqref{e:UH}.
\end{Definition}

\begin{Corollary}\label{c:symp}
The map $\UU$ preserves the symplectic pairing. In particular,
\[
\br{f(-z), g(z)}_{X_-} = \br{\UU f(-z), \UU g(z)}_{X_+}.
\]
\end{Corollary}

\begin{proof}
By the Leray--Hirsch theorem, any class in $\hat H^*_\torus(X_-)$ may be decomposed as a sum of classes of the form $\pi_-^*(\beta) \ch^\torus(\tilde v_-^* \ol C)$ for some $\ol C \in K^0_\torus\bigl(\ol X_-\bigr)$, where the map $\tilde v_-$ appears in the proof of Corollary~\ref{cv}.

We may therefore assume without loss of generality that
$f(z) =\psi_-\bigl( \pi_-^*(\beta_1)\ch^\torus\bigl(\tilde v_-^* \ol C\bigr)\bigr)$ and
$g(z) = \psi_-\bigl(\pi_-^*(\beta_2) \ch^\torus\bigl(\tilde v_-^* \ol D\bigr)\bigr)$.
Using the same chain of equalities from the previous proof, the projection formula, and \eqref{vgdiag}, one computes that
\begin{align*}
 \bigl\langle\psi_-\bigl( \pi_-^*(\beta_1)\ch^\torus\bigl(\tilde v_-^* \ol C\bigr)\bigr)\bigl( {\rm e}^{-\pi {\rm i}} z\bigr),\psi_-\bigl(\pi_-^*(\beta_2) \ch^\torus\bigl(\tilde v_-^* \ol D\bigr)\bigr)(z)\bigr\rangle_{X_-} = K \cdot \chi_z\bigl( \ol C, \ol D\bigr)\end{align*}
and
\begin{align*}
 & \bigl\langle \UU\psi_- \bigl( \pi_-^*(\beta_1)\ch^\torus\bigl(\tilde v_-^* \ol C\bigr)\bigr)\bigl( {\rm e}^{-\pi {\rm i}} z\bigr), \UU\psi_-\bigl(\pi_-^*(\beta_2) \ch^\torus\bigl(\tilde v_-^* \ol D\bigr)\bigr)(z)\bigr\rangle_{X_+}\\&\qquad{} =K \cdot \chi_z\bigl(\ol \FM \bigl(\ol C\bigr), \ol \FM\bigl( \ol D\bigr)\bigr),
 \end{align*}
where
\[
K = \biggl( z^{- \lambda \partial \lambda } \frac{(2 \pi)^{\dim X}}{(2 \pi {\rm i})^{\dim B}} \int_B \bigl((-2\pi {\rm i})^{\deg_0/2} \beta_1\bigr) \cup \bigl((2\pi {\rm i})^{\deg_0/2}\beta_2 \op{Td}_B\bigr) \biggr).
\]
Because $\ol \FM$ is induced from an equivalence of categories,
\[
\chi_z\bigl( \ol C, \ol D\bigr) = \chi_z\bigl(\ol \FM \bigl(\ol C\bigr), \ol \FM\bigl( \ol D\bigr)\bigr).
\]
The claim follows.
\end{proof}

\subsection[I-functions for Grassmann bundles]{$\boldsymbol{I}$-functions for Grassmann bundles}
In this section, we record $I$-functions for $X_\pm$ and the associated abelian quotients. In each case, an $I$-function is obtained by modifying the $J$-function of the base $B$.

\subsubsection[r=1]{$\boldsymbol{r=1}$}
We begin with the case $r=1$,
in which
\begin{align*}
X_{+} = \tot \bigl( \cc O_{\PP(\E)}(-1) \otimes p^*\F^\vee\bigr),\qquad
X_{-} =\tot \bigl( \cc O_{\PP(\F^\vee)}(-1) \otimes p^*\E\bigr).
\end{align*}
In this case $X_{\pm}$ is a toric bundle over $B$. $I$-functions for split toric bundles were computed by Brown in \cite{Brown} in terms of the $J$-function for the base $B$. Let $J^B = \sum_{\beta \in NE(B)_\ZZ} Q^\beta J_\beta$ be the $J$-function for $B$. Define the $H^*_\torus(X_\pm)$-valued function $I_{X_\pm}(Q, q, z)$ to be
\begin{align*}
&I_{X_+}(Q, q_+, z) =\sum_{\beta \in NE(B)_\ZZ} Q^\beta J_\beta \sum_{e \in \ZZ} q_+^{-\y/z+e} M_{\beta, -e} ,\\
&I_{X_-}(Q, q_-, z) =\sum_{\beta \in NE(B)_\ZZ} Q^\beta J_\beta \sum_{d \in \ZZ} q_-^{\y/z+d} M_{\beta, d},
\end{align*}
where
\begin{equation}\label{hmod} M_{\beta, d} =
\prod_{i=1}^n \frac{\prod_{h=-\infty}^0 (\z_i+\y+hz)}{\prod_{h=-\infty}^{d-\beta\cdot \MB_i} (\z_i+\y+hz)}\prod_{k=1}^{n}\frac{\prod_{h=-\infty}^{0}(-\x_k -\y+hz)}{\prod_{h=-\infty}^{-d+\beta\cdot \LB_k}(-\x_k -\y+hz)}.\end{equation}
\begin{Remark}
Despite having the same formal expression under the change of variables $q_+ = q_-^{-1}$, we emphasize that $I_{X_+}$ and $I_{X_-}$ take very different forms due to the fact that they lie in different cohomology rings. In particular, in $I_{X_+}$, the expression $M_{\beta, -e}$ vanishes whenever $e < \op{min}_k(-\beta\cdot \LB_k)$ by \eqref{e:rels}. On the other hand, in $I_{X_-}$ the expression $M_{\beta, d}$ vanishes when $d < \op{min}_i(\beta\cdot \MB_i)$.

\end{Remark}

The following is a special case of Brown's result for toric bundles (see \cite{Brown}).
\begin{Theorem}
The function $I_{X_\pm}(Q, q_\pm, -z)$ lies in the $\torus$-equivariant Lagrangian cone $\cc L_{X_\pm}$.
\end{Theorem}

\subsubsection{The associated abelian quotients}
Now assume $r > 1$. By iterating the result of the previous section, we obtain $I$-functions for~$X_{T, \pm}$, and for all $X_{T, j}$ for $0 \leq j \leq s$.
Let
\begin{align*}
	L_{\beta,\vec{d}}={}& \prod_{j=1}^r \Biggl(\prod_{i=1}^n \frac{\prod_{h=-\infty}^0 (\z_i+\y_j+hz)}{\prod_{h=-\infty}^{d_j-\beta\cdot \MB_i} (\z_i+\y_j+hz)}
	\cdot \prod_{k=1}^{n}\frac{\prod_{h=-\infty}^{0}(-\x_k -\y_j+hz)}{\prod_{h=-\infty}^{-d_j+\beta\cdot \LB_k}(-\x_k -\y_j+hz)}\Biggr).
\end{align*}
Define $I_{X_{T, j}} (Q, d_1, \dots, d_j, e_{j+1}, \dots, e_r, z)$ to be
\begin{equation*} I_{X_{T, j}} =
\sum_{\beta \in NE(B)_\ZZ} Q^\beta J_\beta \sum_{\stackrel{(e_1, \dots, e_j) \in \ZZ^j}{ (d_{j+1}, \dots, d_r) \in \ZZ^{s-j}}} \prod_{i=1}^j q_{i, +}^{-\y_i/z+e_i}\prod_{i=j+1}^r q_{i, -}^{\y_i/z+d_i}
L_{\beta, d_1, \dots, d_j, -e_{j+1}, \dots, -e_r}.
\end{equation*}
Again by \cite{Brown}, we see that $I_{X_{T, j}} (Q, d_1, {\dots}, d_j, e_{j+1}, {\dots}, e_r, -z)$ lies on the Lagrangian cone for~$X_{T, j}$.

\subsubsection{The non-abelian quotients}
The abelian/non-abelian correspondence in Gromov--Witten theory is the principle that one can recover the Gromov--Witten theory of the non-abelian quotient from that of the associated abelian quotient. There is a long history of results in this direction, e.g., \cite{BCK, BCK2, CKS, CLS, Web, Web2}. The modern formulation of the expected relationship is, roughly, that the (genus zero) Gromov--Witten theory for the non-abelian quotient may be obtained from that of the associated abelian quotient by the following steps:
\begin{enumerate}\itemsep=0pt
\item[1)] \emph{twisting} the theory of the abelian quotient by the Euler class of the root bundle $\mathfrak g/\mathfrak t$;
\item[2)] \emph{specializing} those Novikov parameters which are identified via the natural map
\[
H_2(BT) \to H_2(BG).
\]
\end{enumerate}
Let
\begin{gather}
 N_{\beta,\vec{d}}=
	 \prod_{j=1}^r\Biggl(\prod_{i=1}^n \frac{\prod_{h=-\infty}^0 (\z_i+\y_j+hz)}{\prod_{h=-\infty}^{d_j-\beta\cdot \MB_i} (\z_i+\y_j+hz)}
	\cdot \prod_{k=1}^{n}\frac{\prod_{h=-\infty}^{0}(-\x_k -\y_j+hz)}{\prod_{h=-\infty}^{-d_j+\beta\cdot \LB_k}(-\x_k -\y_j+hz)}\Biggr) \nonumber \\
	\hphantom{ N_{\beta,\vec{d}}=\prod_{j=1}^r}{}
 \times \prod_{\substack{i,j=1\\ i\neq j}}^r
	\frac{\prod_{h=-\infty}^{d_i-d_j}(\y_i-\y_j+hz)}{\prod_{h=-\infty}^{0}(\y_i-\y_j+hz)} .\label{Nbeta}
	\end{gather}

Define \smash{$I_{X_{T, +}}^{e_\torus(\mathfrak g /\mathfrak t)}(Q, \vec q_+, z)$} and \smash{$I_{X_{T, -}}^{e_\torus(\mathfrak g /\mathfrak t)}(Q, \vec q_-, z)$} to be
\begin{align*}
&I_{X_{T, +}}^{e_\torus(\mathfrak g /\mathfrak t)}=
\sum_{\beta \in NE(B)_\ZZ} Q^\beta J_\beta \sum_{(e_1, \dots, e_r) \in \ZZ^r} \prod_{i=1}^r q_{i, +}^{-\y_i/z+e_i}N_{\beta, -e_1, \dots , -e_r}, \\
&I_{X_{T, -}}^{e_\torus(\mathfrak g /\mathfrak t)}=
\sum_{\beta \in NE(B)_\ZZ} Q^\beta J_\beta \sum_{(d_1, \dots, d_r) \in \ZZ^r} \prod_{i=1}^r q_{i, -}^{\y_i/z+d_i}N_{\beta, d_1, \dots , d_r}.
\end{align*}
Then, define
\[
I_{X_\pm}(Q, q_\pm, z) = I_{X_{T, \pm}}^{e_\torus(\mathfrak g /\mathfrak t)}(Q, \vec q_\pm, z)|_{q_{1, \pm} = \cdots = q_{s, \pm} = q_\pm}.
\]
The following proposition is a corollary of Coates--Lutz--Shafi in \cite{CLS} on the abelian/non-abelian correspondence. One could also prove it directly by generalizing the work of Oh in \cite{Oh1} on $I$-functions of Grassmann bundles.
\begin{Proposition}
The function $I_{X_\pm}(Q, q_\pm, -z)$ lies on the $\torus$-equivariant Lagrangian cone $\cc L_{X_\pm}$.
\end{Proposition}
\begin{proof}

As described in \cite{CLS}, $I_{X_\pm}$ is the \emph{Givental--Martin} modification of $I_{X_{T, \pm}}$, so the proposition essentially follows from \cite[Theorem 1.11]{CLS}.
More precisely, we apply the argument of~\mbox{\cite[Theorem~5.11]{CLS}}, which provides an abelian/non-abelian correspondence for twisted theories over a flag bundle. Because that theorem is not stated in exactly the generality which we need here, we summarize the argument below and explain the necessary adjustments to the case at hand. We will prove the proposition for $X_+$.

The space $X_+$ is the total space of a vector bundle over the Grassmann bundle
\[
Z := \Gr(r, \E) = \tot\left(\mathcal{H}{\rm om}(\cc O_B \otimes \CC^r, \E)\right) \sslash_+ G.
\]
The associated abelian quotient is the toric bundle given by the $r$-fold fiber product of $\PP(\E)$ over $B$:
\[
Z_T := \tot\left(\mathcal{H}{\rm om}(\cc O_B \otimes \CC^r, \E)\right) \sslash_+ T = \PP(\E) \times_B \cdots \times_B \PP(\E).
\]

Let $I_{Z_T}$ denote the Brown $I$-function for $Z_T$ as in \cite{Brown}, and let
\smash{$I_{Z_T}^{e_\torus^{-1}( (S_T \otimes p^*\F^\vee))}$} denote the corresponding
\smash{$\bigl(\bigl(S_T \otimes p^*\F^\vee\bigr), {\rm e}^{-1}_{\torus}(-)\bigr)$}-twisted $I$-function, where we use $S_T$ to denote the vector bundle on $Z_T$ induced by the right standard representation $R$ from \eqref{e:srep} restricted to $T \subset G$. Then \raisebox{-1pt}{\smash{$I_{Z_T}^{e_\torus^{-1}( (S_T \otimes p^*\F^\vee))}$}} is exactly the $I$-function $I_{X_{T,+}}$ from the previous section.

Following the notation and discussion in \cite[Section~5.2]{CLS} (in particular, Theorem~5.11), we apply a further twist to
\[
I_{X_{T,+}} = I_{Z_T}^{ {\rm e}^{-1}_{\torus}( (S_T \otimes p^*\F^\vee)}
\]
by $(\Phi, \mathbf c')$ where $\Phi \to Z_T$ is the root bundle with respect to $G$ and $\mathbf c'$ is the $\CC^*$-equivariant Euler class
with respect to the $\CC^*$-action scaling the fibers of $\Phi$. Denote the equivariant parameter of this $\CC^*$ action by $\lambda$. Specializing variables to
$q_{1, +} = \cdots = q_{s, +} = q_+,$
projecting via the map on Givental spaces defined in \cite{CLS}
\[
p\colon \ \mathcal{H}_{Z_T}^W \to \mathcal{H}_{Z},
\]
and taking the nonequivariant limit $\lambda \mapsto 0$, we obtain $I_{X_+}(Q, q_+, z)$.

By the proof of \cite[Theorem~5.11]{CLS}, this function lies in the $\bigl(\bigl(S_+ \otimes p^*\F^\vee\bigr), {\rm e}^{-1}_{\torus}\bigr)$-twisted Lagrangian cone of $Z = \Gr(r, \E)$. While the theorem in loc.~cit.\ is only stated for the twisted theory $( L, e_{\CC^*}(-))$, where $L$ is a line bundle and $e_{\CC^*}(-)$ is an equivariant Euler class, by~\mbox{\cite[Remark 3.7]{CLS}} the result holds for other multiplicative classes as well. The further extension to the case that~$L$ is replaced by the sum of line bundles also follows the same argument.

By definition, the $\bigl(\bigl(S_+ \otimes p^*\F^\vee\bigr), {\rm e}^{-1}_{\torus}\bigr)$-twisted Lagrangian cone of $\Gr(r, \E)$ is equal to the $\torus$-equivariant Lagrangian cone of $X_+$.
\end{proof}

\section{Analytic continuation: The abelian case}

In this section, we explicitly compute the analytic continuation of the $I$-functions of the associated abelian quotients $X_{T, \pm}$ and show that they are related by symplectic transformation. This was proven in \cite{CIJ} in the case $B$ is a point.

\subsection{Projective bundles}
Let us consider the special case of $r=1$ first. We will analytically continue the function $I_{X_+}$ to $q_+ = \infty$, and show that the resulting function is related to $I_{X_-}$ by a symplectic transformation.

We begin by writing $I_{X_\pm}$ as $\psi (H_{X_\pm})$ where $\psi$ is as in \eqref{psi}.
By \eqref{e:Euler} in the case $r=1$ we obtain the relative Euler sequence
\[
0 \to \cc O_{X_-} \to S_- \otimes p^*E \oplus S_-^\vee \otimes p^*F^\vee \to TX_{-}/B \to 0
\]
and similarly for $X_+$.
Because the Gamma class is multiplicative, we have that
\begin{align*}
\hat \Gamma_{X_\pm} = \hat \Gamma_{B} \hat \Gamma(TX_\pm/B)
 =\hat \Gamma_{B} \prod_{j=1}^n \Gamma(1 + (- \x_j - \y)) \prod_{i=1}^n \Gamma(1 + (\z_i + \y)).
\end{align*}
Note that the hypergeometric modification \eqref{hmod} is equal to
\begin{equation*}\frac{1}{z^{\beta \cdot c_1(\E \oplus \F^\vee)}}\prod_{j=1}^n\frac{ \Gamma(1 + (- \x_j - \y)/z) \Gamma(1 + (\z_j + \y)/z)}{ \Gamma(1 + (- \x_j - \y)/z + \beta \cdot \LB_j-d) \Gamma(1 + (\z_j + \y)/z - \beta \cdot \MB_j + d)}.
\end{equation*}

Define $H^B(Q ,z)$ to be
\begin{align*} H^B := (2 \pi {\rm i})^{-\deg_0/2}\frac{1}{\hat \Gamma_B} z^{- \rho_{X}} z^{\mu_{X}} J^B
 =
\sum_{\beta \in NE(B)_\ZZ} Q^\beta H_{B, \beta},\end{align*}
where $\rho_{X}$ is the operator on $H^*_\torus(B)$ given by multiplication by
$c_1^\torus\bigl(TB \oplus E \oplus F^\vee\bigr)$ and
\[
\mu_{X}(\phi) = \biggl( \frac{\deg(\phi)}{2} -\frac{ \dim(X_\pm)}{2} \biggr) \phi.
\]

Define the $H^*_\torus\bigl(X_+\bigr)$-valued function
\begin{align*}&H_{X_{+}, \beta}(q_+)\\ \nonumber
&\qquad{}= \sum_{e \in \ZZ} \frac{q_+^{-y/2 \pi {\rm i}+e}}
{\prod_{j=1}^n \Gamma(1 + (- \x_j - \y)/2 \pi {\rm i} + \beta \cdot \LB_j+e) \Gamma(1 + (\z_j + \y)/2 \pi {\rm i} - \beta \cdot \MB_j - e)}.\end{align*}
Similarly, define the $H^*_\torus(X_-)$-valued function
\begin{align}\label{e:HXm} &H_{X_-, \beta}(q_-) \\ \nonumber
& \qquad{}= \sum_{d \in \ZZ} \frac{q_-^{y/2 \pi {\rm i}+d}}
{\prod_{j=1}^n \Gamma(1 + (- \x_j - \y)/2 \pi {\rm i} + \beta \cdot \LB_j-d) \Gamma(1 + (\z_j + \y)/2 \pi {\rm i} - \beta \cdot \MB_j + d)}.\end{align}

Then a simple check shows that
\begin{align}\label{ItoH} I_{X_\pm}(Q, q_\pm, z) &{}= \psi \biggl( \sum_{\beta \in NE(B)_\ZZ} \frac{Q^\beta}{z^{\beta \cdot c_1(\E \oplus \F^\vee)}} H_{B, \beta} H_{X_\pm, \beta} \biggr).
\end{align}

We will fix $\beta$ and relate $H_{X_+, \beta}$ and $H_{X_-, \beta}$ via analytic continuation and symplectic transformation.
Recall Definition~\ref{d:UH0} in the special case $r=1$. In this case, $\UU_H\colon H^*_\torus(X_{-}) \otimes_{R_\torus} S_\torus \to H^*_\torus(X_{+})\otimes_{R_\torus} S_\torus$ is the linear map defined
by sending
$\alpha/e_\torus(N_{B_k|X_-}) \in H^*_\torus(B_k^{-})$ to
\begin{equation}\label{eU1}
\sum_{l = 1}^n {\rm e}^{ (n-1)\cdot ( \x_l-\z_k)/2 }\frac{ \prod_{i\neq k} \sin( (-\z_i+\x_l)/2 {\rm i}) }{\prod_{i \neq k}\sin((-\z_i + \z_k)/2 {\rm i})} \phi_{kl}(\alpha)/e_\torus(N_{B_l|X_+}),\end{equation}
where $\phi_{kl}\colon H^*_\torus(B_k^{-})\otimes_{R_\torus} S_\torus \to H^*_\torus(B_l^{+})\otimes_{R_\torus} S_\torus$ is the canonical isomorphism induced by the projections $\pi_k^{+}\colon B_k^{+} \to B$ and $\pi_l^{-}\colon B_l^{-} \to B$.
By restricting $H_{X_+, \beta}$ to the fixed points of the $\torus$-action and
using Mellin--Barnes integrals, we may analytically continue $H_{X_+, \beta}$ to $q_+ = \infty.$ Because the method of computation is by now somewhat standard (see, e.g., \cite{CIJ}), we relegate the proof to Appendix~\ref{appendix}.
We obtain the following.
\begin{Proposition}\label{p:ac1}
We have the equality
 \begin{equation}\label{analyticcomp}
\UU_H H_{X_{-}, \beta} = \widetilde{H_{X_{+}, \beta}}
\end{equation}
for all $\beta$, where $\widetilde{H_{X_{+}, \beta}}$ denotes the analytic continuation of $H_{X_{+}, \beta}$ along the path described above, under the substitution $q_+^{-1} = q_-$.
\end{Proposition}

Recall Definition~\ref{d:UH}, $\UU = \psi_{+} \circ \UU_H \circ \psi_{-}^{-1}.$ We conclude the following in the case $r=1.$
\begin{Theorem} \label{t:st}
The linear transformation $\UU$ is symplectic, has a well-defined non-equivariant limit, and satisfies:
\[
\UU I_{X_{-}} =\widetilde{I_{X_{+}}},
\]
where $\widetilde{I_{X_{+}}}$ denotes the analytic continuation of $I_{X_{+}}$ along the path $\gamma$ and we set $q_+^{-1} = q_-$.
Furthermore, the following diagram commutes:
\[
\begin{tikzcd}
K^0_\torus(X_{-}) \ar[r, "\FM"] \ar[d, "\Psi_-"] & K^0_\torus(X_{+}) \ar[d, "\Psi_+"] \\
\widetilde{\cc H}_{X_{-}} \ar[r, "\UU"] & \widetilde{\cc H}_{X_{+}}.
\end{tikzcd}
\]
\end{Theorem}

\begin{proof}
The proof of the equality is immediate from Definition~\ref{d:UH}, \eqref{ItoH} and \eqref{analyticcomp}. The fact that~$\UU$ is symplectic and has a well-defined non-equivariant limit is Corollary~\ref{c:symp}. The commuting diagram follows from Proposition~\ref{p:Kcom}.
\end{proof}

\subsection{Products of projective bundles}
We next consider, for general $r >0$, the wall crossing between $I$-functions of the associated abelian quotients $X_{T, \pm}$. While the result follows simply by iterating Theorem~\ref{t:st} $r$ times, the computations in this section will play a role when we consider the wall crossing of the Grassmann~flop.

To pass from $X_{T, -}$ to $X_{T, +}$ via variation of GIT, one must cross $r$ distinct walls, one for each factor in $T = (\CC^*)^r$.
Recall the definition of $X_{T, j}$ from Section~\ref{s:aa}, and set $X_{T, -} = X_{T, 0}$ and $X_{T, +} = X_{T, r}$.
Observe that the wall crossing from $X_{T, j}$ to $X_{T, j+1}$ is a particular example of the wall crossing from the previous section.

In this case, for any $j$ we have the relationship between Gamma functions:
\[
\hat \Gamma(TX_{T, j})=
\hat \Gamma_{B} \hat \Gamma(TX_T/B)
\]
with
\[
\hat \Gamma(TX_T/B) = \prod_{k=1}^r\prod_{j=1}^n\bigl(
\hat \Gamma(1+(- y_k-\x_j))\Gamma(1+(y_k+\z_j))\bigr).
\]

Define $H^B_T(Q ,z)$ to be
\begin{align*} H^B_T := (2 \pi {\rm i})^{-\deg_0/2}\frac{1}{\hat \Gamma_B} z^{- \rho_{X_T}} z^{\mu_{X_T}} J^B
 =
\sum_{\beta \in NE(B)_\ZZ} Q^\beta H_{B, T, \beta},\end{align*}
where $\rho_{X_T}$ is the operator on $H^*_\torus(B)$ given by multiplication by \smash{$c_1^\torus\bigl(TB \oplus \bigl(E \oplus F^\vee\bigr)^{\oplus r}\bigr)$} and
\[
\mu_{X_T}(\phi) = \biggl( \frac{\deg(\phi)}{2} -\frac{ \dim(X_{T, \pm})}{2} \biggr) \phi.
\]

Define the $H^*_\torus(X_{T,+})$-valued, resp.\ $H^*_\torus(X_{T,-})$-valued, functions
\begin{align*}
&H_{X_{T, +}, \beta} =\\ &
\sum_{(e_1, \dots, e_r) \in \ZZ^r} \prod_{k=1}^r\prod_{j=1}^n \!\frac{q_{k, +}^{(-\y_k/2 \pi {\rm i} + e_k)}}
 {\Gamma(1+(-\y_k-\x_j)/2 \pi {\rm i}+\beta\cdot \LB_j+e_k)\Gamma(1+(\y_k+\z_j)/2 \pi {\rm i}-\beta\cdot \MB_j-e_k)},
\\
&H_{X_{T, -}, \beta} = \\ & \sum_{(d_1, \dots, d_r) \in \ZZ^r} \prod_{k=1}^r \prod_{j=1}^n\!\frac{q_{k, -}^{(\y_k/2 \pi {\rm i} + d_k)} }{\Gamma(1+(-\y_k-\x_j)/2 \pi {\rm i}+\beta\cdot \LB_j - d_k)\Gamma(1+(\y_k+\z_j)/2 \pi {\rm i}-\beta\cdot \MB_j + d_k)}.
\end{align*}

Then as in the case of $r=1$, we have
\begin{align*}
I_{X_{T, \pm}}(Q, q_\pm, z) &{}= \psi_{X_{T,\pm}} \biggl( \sum_{\beta \in NE(B)_\ZZ} \frac{Q^\beta}{z^{\beta \cdot c_1^\torus(\E \oplus \F^\vee)^{\oplus r}}} H_{B, T, \beta} H_{X_{T, \pm}, \beta} \biggr).
\end{align*}

Recall that the $\torus$-fixed loci of $X_{T, \pm}$ are indexed by $F_\pm$.

Let $\UU_H^T\colon H^*_\torus(X_{T,-})\otimes_{R_\torus} S_\torus \to H^*_\torus(X_{T,+})\otimes_{R_\torus} S_\torus$ be the map defined by sending
$\alpha/e_\torus\bigl(N_{f^-}\bigr) \in H^*_\torus\bigl(B_{f^-}\bigr)$ to
\[
\sum_{f^+ \in F_+} \prod_{i=1}^r {\rm e}^{ (n-1)\cdot ( \x_{f^+_i}- \z_{f^-_i})/2 }\frac{ \prod_{j\neq {f^-_i}} \sin\bigl( \bigl(-\z_j+\x_{f^+_i}\bigr)/2 {\rm i}\bigr) }{\prod_{j \neq {f^-_i}}\sin\bigl(\bigl(-\z_j + \z_{f^-_i}\bigr)/2 {\rm i}\bigr)} \frac{\phi_{f^+, f^-}(\alpha)}{e_\torus(N_{f^+})},
\]
where $\phi_{f^+, f^-}\colon H^*_\torus\bigl(B_{f^-}\bigr)\otimes_{R_\torus} S_\torus \to H^*_\torus\bigl(B_{f^+}\bigr)\otimes_{R_\torus} S_\torus$ is the canonical isomorphism induced by the projections $\pi_{f^+}\colon B_{f^+} \to B$ and $\pi_{f^-}\colon B_{f^-} \to B$, and $N_{f^\pm}$ is the normal bundle $N_{B_{f^\pm}|X_\pm}$.

A straightforward generalization of the analytic continuation in Appendix~\ref{appendix} (essentially repeating it $r$ times, one for each $q_i$ variable) yields the following:
\[
\UU_H^T H_{X_{T,-}, \beta} = \widetilde{H_{X_{T,+}, \beta}}
\]
for all $\beta$, after identifying $q_{k, +}^{-1} = q_{k, -}$. Here the path of analytic continuation is a concatenation
$\gamma_1^T \star \cdots \star \gamma_r^T$,
where $\gamma_k^T$ is the path that keeps $q_{i, +}$ constant for $i \neq k$, and in the $ \log(q_{k, +})$ plane follows the path $\gamma$ of analytic continuation from the previous section.

Define
$\UU^T = \psi_{X_{T,+}} \circ \UU_H^T \circ \psi_{X_{T,-}}^{-1}$
in analogy to the previous section. By Theorem~\ref{t:st}, we obtain the following.
\begin{Theorem} 
The linear transformation $\UU^T$ is symplectic, has a well-defined non-equivariant limit, and satisfies
\[
\UU^T I_{X_{T,-}} =\widetilde{I_{X_{T,+}}}.
\]
\end{Theorem}

\section{Grassmann flops}
In this section, we compute the analytic continuation and symplectic transformation for the $I$-functions $X_+$ and $X_-$ appearing in the Grassmann flop.

\subsection{Analytic continuation}
The computation proceeds by first computing the analytic continuation for the \emph{twisted} $I$-functions of the abelian quotients, similar to the previous section, and then deforming the path to the locus $q_{1, +} = \cdots = q_{r, +}$.

Let $S_-$ denote the rank-$r$ tautological bundle on $X_-$. By \eqref{e:Euler}, we have the short exact sequence
\[0 \to S_- \otimes S_-^\vee \to S_- \otimes p^*E \oplus S_-^\vee \otimes p^*F^\vee \to TX_-/B \to 0
\]
and similarly for $X_+$. We view $H^*_\torus(X_\pm)$ as a subspace of $\bigl(H^*_\torus( V^{s}(G)/T)\bigr)$ via Proposition~\ref{p:pull}. After pulling back to $V^s(G)/T$, we may
 write the Gamma class $\hat \Gamma_{X_\pm}$ as
$\hat \Gamma_{B} \hat \Gamma(TX_\pm/B)$, where in this case
\[
\hat \Gamma(TX_\pm/B) = \frac{\prod_{k=1}^r\prod_{j=1}^n\bigl(
\hat \Gamma(1+(- y_k-\x_j))\Gamma(1+(y_k+\z_j))\bigr)}{\hat \Gamma(\mathfrak g/\mathfrak t)}
\]
with
\[
\hat \Gamma(\mathfrak g/\mathfrak t) = \prod_{i \neq k} \hat \Gamma(1 + y_k - y_i).
\]
Recall that the $I$-functions for $X_\pm$ are specializations at $q_{1, \pm} = \cdots = q_{r, \pm}$ of \smash{$I_{X_{T, \pm}}^{e_\torus(\mathfrak g /\mathfrak t)}$} where
\begin{align*}
&I_{X_{T, +}}^{e_\torus(\mathfrak g /\mathfrak t)}=
\sum_{\beta \in NE(B)_\ZZ} Q^\beta J_\beta \sum_{(e_1, \dots, e_r) \in \ZZ^r} \prod_{i=1}^r q_{i, +}^{-\y_i/z+e_i}N_{\beta, -e_1, \dots , -e_r} ,\\
 &I_{X_{T, -}}^{e_\torus(\mathfrak g /\mathfrak t)}=
\sum_{\beta \in NE(B)_\ZZ} Q^\beta J_\beta \sum_{(d_1, \dots, d_r) \in \ZZ^r} \prod_{i=1}^r q_{i, -}^{\y_i/z+d_i}N_{\beta, d_1, \dots , d_r}
\end{align*}
with $N_{\beta,\vec{d}}$ defined in \eqref{Nbeta}.

The hypergeometric factor $N_{\beta,\vec{d}}$ can be rewritten as \begin{align*}
 &\frac{1}{z^{\beta \cdot c_1^\torus(\E \oplus \F^\vee)^{\oplus r}}}\prod_{ i\neq k}\frac{ \Gamma(1 + (\y_k - \y_i)/z + d_k - d_i) }{ \Gamma(1 + (\y_k - \y_i)/z)} \\
 &\qquad{}\times\prod_{k=1}^r\prod_{j=1}^n \frac{
\hat \Gamma(1+(- \y_k-\x_j)/z)\Gamma(1+(\y_k+\z_j)/z)}{\Gamma(1+(-\y_k-\x_j)/z+\beta\cdot \LB_j - d_k)\Gamma(1+(\y_k+\z_j)/z-\beta\cdot \MB_j + d_k)}.
\end{align*}

Generalizing the definition of $H^B(Q ,z)$ to the case of $r>1$, we define
\begin{align*} H^B := (2 \pi {\rm i})^{-\deg_0/2}\frac{1}{\hat \Gamma_B} z^{- \rho_{X}} z^{\mu_{X}} J^B
 =
\sum_{\beta \in NE(B)_\ZZ} Q^\beta H_{B, \beta},\end{align*}
where $\rho_{X}$ is the operator on $H^*_\torus(B)$ given by multiplication by $c_1^\torus\bigl(TB \oplus \bigl(E \oplus F^\vee\bigr)^{\oplus r}\bigr)$ and $\mu_X$ acts on a homogeneous element $\phi$ by
\[
\mu_{X}(\phi) = \biggl( \frac{\deg(\phi)}{2} -\frac{ \dim(X_{\pm})}{2} \biggr) \phi.
\]

Denote by $H_{X_+, \beta}$ the function
\begin{gather*}
 \sum_{(e_1, \dots, e_r) \in \ZZ^r} \prod_{k=1}^r\! \Biggl(\!q_{k, +}^{(-\y_k/2\pi {\rm i} + e_k)}
\!\prod_{i < k}\! \Gamma(1 + (\y_k - \y_i)/2 \pi {\rm i} - e_k + e_i) \Gamma(1 - (\y_k - \y_i)/2 \pi {\rm i} + e_k - e_i)\! \\
\quad\quad{}\times \prod_{j=1}^n\Gamma(1+(-\y_k-\x_j)/2 \pi {\rm i}+\beta\cdot \LB_j+e_k)^{-1}\Gamma(1+(\y_k+\z_j)/2 \pi {\rm i}-\beta\cdot \MB_j-e_k)^{-1}\!\Biggr)
\\
\quad{} = \frac{\pi^{r(r-1)/2}}{\prod_{i<k\leq r} \sin((\y_k - \y_i)/2{\rm i})} \sum_{(e_1, \dots, e_r) \in \ZZ^r} (-1)^{(r-1) \sum_{k=1}^r e_k}
 \\
\qquad{} \times \prod_{k=1}^r \Biggl(q_{k, +}^{(-\y_k/2\pi {\rm i} + e_k)}
\prod_{i < k} ((\y_k - \y_i)/2 \pi {\rm i} - e_k + e_i)
 \\
\qquad{}\times \prod_{j=1}^n\Gamma(1+(-\y_k-\x_j)/2 \pi {\rm i}+\beta\cdot \LB_j+e_k)^{-1}\Gamma(1+(\y_k+\z_j)/2 \pi {\rm i}-\beta\cdot \MB_j-e_k)^{-1}\Biggr).
\end{gather*}
As before, we have that
\begin{align}\label{ItoH3} I_{X_+}(Q, q_+, z) &{}= \psi_+ \left( \sum_{\beta \in NE(B)_\ZZ} \frac{Q^\beta}{z^{\beta \cdot c_1^\torus(\E \oplus \F^\vee)}} H_{B, \beta} H_{X_+, \beta} \right)\vline_{q_{1, +} = \cdots =q_{r, +} = q_+}.
\end{align}
An analogous statement holds for $X_-$.

Restrict $H_{X_+, \beta}$ to a fixed locus $B_\delta$, where $\delta \in D_+$ is a size $r$ subset of $\{1, \dots, n\}$. Then
\begin{equation}\label{HBP}
 H^{+}_{\beta, \delta} := H_{X_+, \beta}|_{B_\delta} =
\frac{\pi^{r(r-1)/2}}{\prod_{i<k\leq r} \sin((\x_{\delta_i} - \x_{\delta_k})/2{\rm i})} K^{+}_{\beta, \delta},
\end{equation}
where $K^{+}_{\beta, \delta}$ is
\begin{align*}
&\sum_{e_1 \geq - \beta \cdot \LB_{\delta_1}} \cdots \sum _{e_r \geq - \beta \cdot \LB_{\delta_r}} (-1)^{(r-1) \sum_{k=1}^r e_k}\prod_{k=1}^r \Biggl(q_{k, +}^{( \x_{\delta_k}/2\pi {\rm i} + e_k)}
\prod_{i < k} ((- \x_{\delta_k} + \x_{\delta_i})/2 \pi {\rm i} - e_k + e_i)
 \\
& \quad{}\times \prod_{j=1}^n\Gamma (1+ (\x_{\delta_k}-\x_j )/2 \pi {\rm i}+\beta\cdot \LB_j+e_k )^{-1}\Gamma (1+ (-\x_{\delta_k}+\z_j )/2 \pi {\rm i}-\beta\cdot \MB_j-e_k )^{-1}\Biggr).
\end{align*}

For any $f \in F_-$, define
 $K^{-}_{\beta, f}$ to be
\begin{align*}
&\sum_{d_1 \geq \beta \cdot \MB_{f_1}} \cdots \sum _{d_r \geq \beta \cdot \MB_{f_r}} (-1)^{(r-1) \sum_{k=1}^r d_k}\prod_{k=1}^r \Biggl(q_{k, -}^{(- \z_{f_k}/2\pi {\rm i} + d_k)}
\prod_{i < k} ( (- \z_{f_k} + \z_{f_i} )/2 \pi {\rm i} + d_k - d_i )
 \\
&\quad{}\times \prod_{j=1}^n\Gamma (1+ ( \z_{f_k}-\x_j )/2 \pi {\rm i}+\beta\cdot \LB_j-d_k )^{-1}\Gamma (1+ (-\z_{f_k}+\z_j )/2 \pi {\rm i}-\beta\cdot \MB_j+d_k )^{-1}\Biggr).
\end{align*}
Then for $\delta^- \in D_-$, after choosing an ordering on the elements and viewing $\delta^-$ as a function, we~have
\begin{equation}\label{HBM}
 H^{-}_{\beta, \delta^-} := H_{X_-, \beta}|_{B_{\delta^-}} =
\frac{\pi^{r(r-1)/2}}{\prod_{i<k\leq r} \sin\bigl(\bigl(\z_{\delta^-_i} - \z_{\delta^-_k}\bigr)/2{\rm i}\bigr)} K^{-}_{\beta, \delta^-}.
\end{equation}

The computation of the analytic continuation of $K_{\beta,\delta}^{+}$ is almost identical to that of $H_{\beta, T, f^+}$ from the previous section.
In this case the path of analytic continuation will be given by $\gamma_1 \star \cdots \star \gamma_r$ where
$\gamma_k$
is the path that keeps $q_{i, +}$ constant for $i \neq k$ and in the
 $\log(q_{k, +})$ plane, will go through $\mathfrak{Im}(\log(q_{k, +})) = (n-r)\pi {\rm i}$ when $\mathfrak{Re}(\log(q_{k,+}))$ is close to 0. See Appendix~\ref{appendix} for details.

Let
\[
C^K_{f^-, f^+} = \prod_{i=1}^r {\rm e}^{ (n-r)\cdot ( \x_{f^+_i}-\z_{f^-_i})/2 }\frac{ \prod_{j\neq {f^-_i}} \sin\bigl( \bigl(-\z_j+\x_{f^+_i}\bigr)/2 {\rm i}\bigr) }{\prod_{j \neq {f^-_i}}\sin\bigl(\bigl(-\z_j + \z_{f^-_i}\bigr)/2 {\rm i}\bigr)}.
\]
Generalizing the computation from Appendix~\ref{appendix},
we see that the analytic continuation of $K^{+}_{\beta, \delta^+}$ along the path
$\gamma_1 \star \cdots \star \gamma_r$ is given by
\begin{equation}\label{ancon}\sum_{f^- \in F_-} C^K_{f^-, \delta^+} K^{-}_{\beta, f^-},\end{equation}
after identifying $q_{k, +}^{-1} = q_{k, -}$. We emphasize that here the sum is over all functions $f^-\colon \{1, \dots, r\} \allowbreak \to \{1, \dots, n\}$.

Then by \eqref{HBP} and \eqref{HBM}, if we define
\[
C^H_{f^-, \delta^+} = C^K_{f^-, \delta^+} \prod_{i<k\leq r} \frac{\sin\bigl(\bigl(\z_{f^-_i} - \z_{f^-_k}\bigr)/2{\rm i}\bigr)}{ \sin\bigl(\bigl(\x_{\delta^+_i} - \x_{\delta^+_k}\bigr)/2{\rm i}\bigr)},
\]
the analytic continuation of $H^{+}_{\beta, \delta^+}$ is given by
\[
\widetilde{H^{+}_{\beta, \delta^+}} =\sum_{f^- \in F_-} C^H_{f^-, \delta^+} H^{-}_{\beta, f^-}.
\]
Note that the numerator of $C^H_{f^-, \delta^+}$ vanishes unless $f^-_i \neq f^-_k$ for $i \neq k$. Therefore, we can replace the summation above to be only over all \emph{injective} functions $\delta^- \in {\rm In}_-$:
\begin{equation}\label{e: CH} \widetilde{H^{+}_{\beta, \delta^+}} = \sum_{\delta^- \in {\rm In}_-} C^H_{\delta^-, \delta^+} H^{-}_{\beta, \delta^-}.\end{equation}

In order to connect the above computation to the actual $I$-functions for $X_+$ and $X_-$ via~\eqref{ItoH3}, we must homotope the path $\gamma_1 \star \cdots \star \gamma_r$ of analytic continuation to a new path which lies entirely within the diagonal $q_{1, +} = \cdots =q_{r, +}$.
\begin{Proposition}\label{p:ancon}
The path of analytic continuation $\gamma_1 \star \cdots \star \gamma_r$ used in \eqref{e: CH} is homotopic to a path $\hat \gamma$ contained entirely in the locus $q_{1, +} = \cdots = q_{r, +}$, via a homotopy which is
\begin{enumerate}\itemsep=0pt
\item[$1)$] contained within the domain for which ${H^{+}_{\beta, \delta^+}}$ is analytic, and
\item[$2)$] independent of $\beta$ and $\delta^+$.
\end{enumerate}
After specializing the Novikov parameters $q_{1, \pm} = \cdots = q_{r, \pm} = q_\pm$, we can therefore analytically continue $H^{+}_{\beta, \delta^+}$ from $q_+ = 0$ to $q_+ = \infty$. Upon setting $q_+^{-1} = q_-$, we obtain the following:
\begin{equation}\label{e: CHspecial} \widetilde{H^{+}}_{\beta, \delta^+}|_{q_{1, +} = \cdots =q_{r, +} = q_+} = \sum_{\delta^- \in {\rm In}_-} C^H_{\delta^-, \delta^+} H^{-}_{\beta, \delta^-}|_{q_{1, -} = \cdots = q_{r, -} = q_-}.\end{equation}
\end{Proposition}
\begin{proof}

First, consider the function $f_{\beta, \delta, k}(q_{k, +})$ defined as
\begin{gather*}
\sum_{e_k \geq - \beta \cdot \LB_{\delta_k}}
\frac{q_{k, +}^{( \x_{\delta_k}/2\pi {\rm i} + e_k)} (-1)^{(r-1)e_k}}{
\prod_{j=1}^n\Gamma(1+(\x_{\delta_k}-\x_j)/2 \pi {\rm i}+\beta\cdot \LB_j+e_k)\Gamma(1+(-\x_{\delta_k}+\z_j)/2 \pi {\rm i}-\beta\cdot \MB_j-e_k)} \\
\qquad{}=\sum_{e_k \geq - \beta \cdot \LB_{\delta_k}}
\left( \frac{q_{k, +}^{( \x_{\delta_k}/2\pi {\rm i} + e_k)}(-1)^{(n+r-1)e_k + \sum_{j=1}^n \beta \cdot M_j }\prod_{j=1}^n \sin((\x_{\delta_k}-\z_j)/2{\rm i})
 }{ \pi^n
} \right. \\
\qquad\quad\hphantom{=\sum_{e_k \geq - \beta \cdot \LB_{\delta_k}}}{} \left.\times \prod_{j=1}^n\frac{\Gamma((\x_{\delta_k}-\z_j)/2 \pi {\rm i}+\beta\cdot \MB_j+e_k)}{\Gamma(1+(\x_{\delta_k}-\x_j)/2 \pi {\rm i}+\beta\cdot \LB_j+e_k)}\right).
\end{gather*}

This satisfies the differential equation
\[
\left[\prod_{j=1}^n (\partial_k -x_j/2 \pi {\rm i} + \beta \cdot L_j) - q_{k, +}(-1)^{n-r+1}\prod_{j=1}^n (\partial_k -z_j/2 \pi {\rm i} + \beta \cdot M_j) \right] \Phi(q_{k,+}) = 0,
\]
where \smash{$\partial_k = q_{k,+} \frac{\partial}{\partial q_{k,+}}$}. This differential equation has singularities at $0$, $\infty,$ and $(-1)^{n-r+1}$. Follo\-wing \eqref{e:cint}, our path of analytic continuation in the $\log(q_{k,+})$ plane goes through $\mathfrak{Im}(\log(q_{k,+})) = (n-r)\pi {\rm i}$ when $\mathfrak{Re}(\log(q_{k,+}))$ is close to 0.
Consider the product
$
\prod_{k=1}^r   f_{\beta, \delta, k}(q_{k, +})$
viewed as an analytic function of the variables $\log(q_{1, +}),\allowbreak \dots, \log(q_{r, +})$. We see that it may be analytically continued everywhere except along the poles:
\begin{equation} \label{e:poles} \log(q_{k,+}) = (n-r+1)\pi {\rm i} + 2 \pi {\rm i} \ZZ.\end{equation}

Define the differential operator
\[
\hat \Delta = \prod_{k=1}^r \biggl( \prod_{i< k} (-\partial_k + \partial_i) \biggr).
\]
Then
\[
K^+_{\beta, \delta} = \hat \Delta \Biggl( \prod_{k=1}^r f_{\beta, \delta, k}(q_{k, +}) \Biggr),
\]
and therefore $ K^+_{\beta, \delta}$ has the same poles at \eqref{e:poles}.

 \begin{figure}[t]\centering
\tikzset{every picture/.style={line width=0.75pt}} 

\begin{tikzpicture}[x=0.75pt,y=0.75pt,yscale=-1,xscale=1]

\draw [line width=3] [dash pattern={on 3.75pt off 18.75pt}] (332.54,10) -- (332.89,143.63) ;
\draw (278.42,89.24) -- (388.73,89.24) ;
\draw (121.24,122.52) .. controls (247.4,123.35) and (227.41,90.9) .. (278.42,89.24) ;
\draw (542.47,123.35) .. controls (420.44,122.52) and (444.57,88.41) .. (388.73,89.24) ;
\draw (138.13,115.2) -- (149.5,122.42) -- (138.13,129.65) ;
\draw (308.58,81.75) -- (319.95,88.98) -- (308.58,96.21) ;
\draw (524.2,115.86) -- (535.57,123.09) -- (524.2,130.32) ;

\draw (335.65,67.33) node [anchor=north west][inner sep=0.75pt] [font=\footnotesize] [align=left] {$\displaystyle ( n-r+1) \pi {\rm i}$};
\draw (335.65,98.12) node [anchor=north west][inner sep=0.75pt] [font=\footnotesize] [align=left] {$\displaystyle {\textstyle ( n-r-1) \pi {\rm i}}$};
\draw (122.46,97.78) node [anchor=north west][inner sep=0.75pt] [font=\footnotesize] [align=left] {$\displaystyle \arg \ =\ 0$};
\draw (491.04,97.95) node [anchor=north west][inner sep=0.75pt] [font=\footnotesize] [align=left] {$\displaystyle \arg \ =\ 0$};
\draw (249.62,65.58) node [anchor=north west][inner sep=0.75pt] [font=\normalsize] [align=left] {$\displaystyle \gamma $};
\end{tikzpicture}
 \caption{}\label{f1}
 \end{figure}

Let $\gamma_1$ denote the
path that keeps $\log(q_{2, +}), \dots, \log(q_{r, +})$ constant (equal to $-\infty$) and follows~$\gamma$ from Figure~\ref{f1} in $\log(q_{1, +})$. Let $\gamma_2$ denote the path that keeps $\log(q_{1, +})$ constant (equal to $\infty$), keeps $\log(q_{3, +}), \dots, \log(q_{r, +})$ constant (equal to $-\infty$) and follows $\gamma$ in $\log(q_{2, +})$. Define $\gamma_k$ for $k \leq r$ similarly. The analytic continuation for $K^+_{\beta, \delta}$ computed in \eqref{ancon} involved the concatenated path
$\gamma_1 \star \cdots \star \gamma_r$.
We must show that this composition of paths is homotopy equivalent to a path $\hat \gamma$ which follows~$\gamma$ in the plane $\log(q_{1, +})= \cdots = \log(q_{r, +})$ via a homotopy which is contained in
\[
\CC^r \setminus \{\log(q_{k,+}) = (n-r+1)\pi {\rm i} + 2 \pi {\rm i} \ZZ\}_{k=1}^r.
\]

We verify this in the case $r=2$. An induction argument proves the claim for general $r$.
First, define a path $\widetilde \gamma_1$ which agrees with $\gamma_1$ except that
\begin{itemize}\itemsep=0pt
\item for all points on the path such that $\mathfrak{Re}(\log(q_{1,+}))\geq 0$ we let $\mathfrak{Im}(\log(q_{1,+})) = (n-r)\pi {\rm i}$;
\item while the value of $\mathfrak{Re}(\log(q_{2,+}))$ remains constant, the value of $\mathfrak{Im}(\log(q_{2,+}))$ is defined to be equal to $\mathfrak{Im}(\log(q_{1,+}))$ for all points on the path.
\end{itemize}

Similarly, define a path $\widetilde \gamma_2$ which agrees with $\gamma_2$ except that
\begin{itemize}\itemsep=0pt
\item for all points on the path such that $\mathfrak{Re}(\log(q_{2,+}))\leq 0$ we let $\mathfrak{Im}(\log(q_{2,+})) = (n-r)\pi {\rm i}$;
 \item while the value of $\mathfrak{Re}(\log(q_{1,+}))$ remains constant, the value of $\mathfrak{Im}(\log(q_{1,+}))$ is defined to be equal to $\mathfrak{Im}(\log(q_{2,+}))$ for all points on the path.
\end{itemize} See Figure~\ref{fig2}.

\begin{figure}\centering
\tikzset{every picture/.style={line width=0.75pt}} 

\begin{tikzpicture}[x=0.75pt,y=0.75pt,yscale=-1,xscale=1]

\draw [line width=3] [dash pattern={on 3.75pt off 18.75pt}] (178.59,12) -- (177.5,120.63) ;
\draw (152.88,90.65) -- (297.5,90.63) ;
\draw (78.24,123.67) .. controls (138.15,124.5) and (128.66,92.3) .. (152.88,90.65) ;
\draw (86.26,116.41) -- (91.66,123.58) -- (86.26,130.76) ;
\draw (167.21,83.22) -- (172.61,90.39) -- (167.21,97.56) ;
\draw (271.6,83.07) -- (277,90.24) -- (271.6,97.42) ;
\draw [line width=3] [dash pattern={on 3.75pt off 18.75pt}] (461.46,10.37) -- (460.5,122.63) ;
\draw (341.5,89.63) -- (490.96,89.62) ;
\draw (571.69,123.72) .. controls (507.61,122.89) and (520.28,88.78) .. (490.96,89.62) ;
\draw (354.37,82.57) -- (360.34,89.8) -- (354.37,97.03) ;
\draw (448.87,82.13) -- (454.85,89.36) -- (448.87,96.58) ;
\draw (562.1,116.24) -- (568.07,123.46) -- (562.1,130.69) ;
\draw (60.5,9.63) -- (318.5,9.63) -- (318.5,135.63) -- (60.5,135.63) -- cycle ;
\draw (337.5,8.63) -- (586.5,8.63) -- (586.5,135.63) -- (337.5,135.63) -- cycle ;

\draw (64.9,99.07) node [anchor=north west][inner sep=0.75pt] [font=\footnotesize] [align=left] {$\displaystyle \arg \ =\ 0$};
\draw (217.94,67.24) node [anchor=north west][inner sep=0.75pt] [font=\footnotesize] [align=left] {$\displaystyle \arg \ =\ ( n-r) \pi {\rm i}$};
\draw (131.22,57.09) node [anchor=north west][inner sep=0.75pt] [font=\normalsize] [align=left] {$\displaystyle \tilde{\gamma }_{1}$};
\draw (338.56,98.15) node [anchor=north west][inner sep=0.75pt] [font=\footnotesize] [align=left] {$\displaystyle \arg \ =\ ( n-r) \pi {\rm i}$};
\draw (532.1,98.32) node [anchor=north west][inner sep=0.75pt] [font=\footnotesize] [align=left] {$\displaystyle \arg \ =\ 0$};
\draw (411.22,53.09) node [anchor=north west][inner sep=0.75pt] [font=\normalsize] [align=left] {$\displaystyle \tilde{\gamma }_{2}$};
\draw (147,138) node [anchor=north west][inner sep=0.75pt] [align=left] {$\displaystyle \log( q_{1,+})$};
\draw (428,138) node [anchor=north west][inner sep=0.75pt] [align=left] {$\displaystyle \log( q_{2,+})$};

\end{tikzpicture}
\caption{The projections of the paths $\widetilde \gamma_1$ and $\widetilde \gamma_2$ to the $\log(q_{1,+})$ and $\log(q_{2,+})$ planes, respectively.}\label{fig2}
\end{figure}
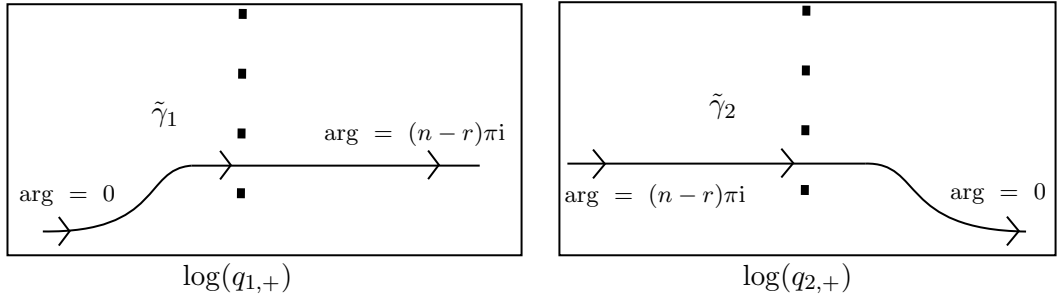

It is clear that
$\widetilde \gamma_1 \star \widetilde \gamma_2$ is homotopic to $ \gamma_1 \star \gamma_2$ in $\CC^r \setminus \{\log(q_{k,+}) = (n-r+1)\pi {\rm i} + 2 \pi {\rm i} \ZZ\}_{k=1}^r.$
By~the~construction, the imaginary parts
$\mathfrak{Im}(\log(q_{1,+}))$ and $\mathfrak{Im}(\log(q_{2,+}))$ are equal for all points on the path $\widetilde \gamma_1 \star \widetilde \gamma_2$.

\begin{figure}\centering
\tikzset{every picture/.style={line width=0.75pt}} 

\begin{tikzpicture}[x=0.75pt,y=0.75pt,yscale=-1,xscale=1]

\draw (269,19.44) -- (391.5,19.44) -- (391.5,136.5) -- (269,136.5) -- cycle ;
\draw [color={rgb, 255:red, 144; green, 19; blue, 254 } ,draw opacity=1 ][line width=1.5] (269,136.5) -- (391.5,136.5) ;
\draw [color={rgb, 255:red, 144; green, 19; blue, 254 } ,draw opacity=1 ][line width=1.5] (391.5,19.44) -- (391.5,136.5) ;
\draw [color={rgb, 255:red, 0; green, 0; blue, 255 } ,draw opacity=1 ][line width=1.5] (391.5,19.44) -- (269,136.5) ;
\draw [color={rgb, 255:red, 144; green, 19; blue, 254 } ,draw opacity=1 ][line width=1.5] (318.5,128.44) -- (332,136.22) -- (318.5,144) ;
\draw [color={rgb, 255:red, 144; green, 19; blue, 254 } ,draw opacity=1 ][line width=1.5] (383.64,87.16) -- (391.08,73.47) -- (399.19,86.78) ;
\draw [color={rgb, 255:red, 0; green, 1; blue, 255 } ,draw opacity=1 ][line width=1.5] (319.99,77.3) -- (335.11,73.53) -- (330.79,88.51) ;

\draw (283,148) node [anchor=north west][inner sep=0.75pt] [align=left] {$\displaystyle \mathfrak{Re}_\torus(\log( q_{1,+}))$};
\draw (403,73) node [anchor=north west][inner sep=0.75pt] [align=left] {$\displaystyle \mathfrak{Re}_\torus(\log( q_{2,+}))$};
\draw (191,128) node [anchor=north west][inner sep=0.75pt] [align=left] {$\displaystyle ( -\infty ,-\infty )$};
\draw (398,11) node [anchor=north west][inner sep=0.75pt] [align=left] {$\displaystyle ( \infty ,\infty )$};

\end{tikzpicture}
\caption{The (real part of the) original path of analytic continuation $\gamma_1 \star \gamma_2 $ (in purple) and the new path $\hat \gamma$ (in blue).}\label{f3}
\end{figure}
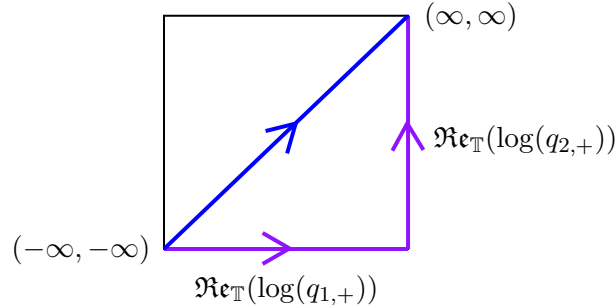

Finally, we homotope $\widetilde \gamma_1 \star \widetilde \gamma_2$ further by modifying the real part of the path so that we have $\mathfrak{Re}(\log(q_{1,+})) = \mathfrak{Re}(\log(q_{2,+}))$ for all points on the path (see Figure~\ref{f3}). This second homotopy is also contained entirely in $\CC^2 \setminus \{\log(q_{k,+}) = (n-2+1)\pi {\rm i} + 2 \pi {\rm i} \ZZ\}_{k=1}^2.$

The resulting path $\hat \gamma$ from
$(\log(q_{1,+}), \log(q_{2,+})) = (-\infty, -\infty)$ to $(\infty, \infty)$ lies entirely in the plane $\log(q_{1,+})= \log(q_{2,+})$ as desired.
\end{proof}

\subsection{Simplification}

After specializing the Novikov parameters, we observe that
$H^{{\pm}}_{\beta, \delta}|_{q_{1, \pm} = \cdots =q_{r, \pm}} = H^{{\pm}}_{\beta, \delta '}|_{q_{1, \pm} = \cdots =q_{r, \pm}}$
whenever $\delta$ and $\delta '$ are equal up to a permutation. This allows for a dramatic simplification to the formula \eqref{e: CHspecial}. The simplification relies on the following identity.
\begin{Lemma}\label{l:sym}
We have the following identity of $($bi$)$-anti-symmetric functions:
\begin{equation}\label{e:sym}
\prod_{i=1}^r \prod_{l<i} ((\z_i - \z_l)(\x_l - \x_i) ) = \sum_{\sigma \in S_r} \op{sgn}(\sigma) \prod_{l=1}^r \prod_{j \neq \sigma(l)} (\x_j - \z_l).
\end{equation}
\end{Lemma}

\begin{proof}
Denote the left- and right-hand sides of \eqref{e:sym} by $L(\underline x, \underline z)$ and $R(\underline x, \underline z)$, respectively.
We first claim that $R(\underline x, \underline z)$ is divisible by $(z_i - z_l)$ for $i \neq l$. Consider the substitution $z_i = z_l = z$.
Group the elements $\sigma \in S_r$ in pairs by grouping $\sigma$ with $\tilde \sigma = \sigma \circ \tau_{il}$ where $\tau_{il}$ is the transposition $(il)$. The elements $\tilde \sigma$ and $\sigma$ have opposite signs. Since $\sigma(l) = \tilde \sigma(i)$ and $\sigma(i) = \tilde \sigma(l)$, we see that
\[
\Biggl(\op{sgn}(\sigma) \prod_{l=1}^r \prod_{j \neq \sigma(l)} (\x_j - \z_l) + \op{sgn}(\tilde \sigma) \prod_{l=1}^r \prod_{j \neq \tilde \sigma(l)} (\x_j - \z_l)\Biggr)\bigg|_{z_i = z_l = z} = 0,
\]
and the claim follows.

On the other hand, $R(\underline x, \underline z)$ may also be written as
\[
\sum_{\sigma \in S_r} \op{sgn}(\sigma) \prod_{j=1}^r \prod_{l \neq \sigma(j)} (\x_j - \z_l),
\]
from which one can apply the same argument to conclude that $R(\underline x, \underline z)$ is divisible by $(x_i - x_l)$ for $i \neq l$.

Because $R(\underline x, \underline z)$ is a polynomial of degree $r(r-1)$, the above divisibility claims imply that it is equal to a constant times $L(\underline x, \underline z)$. We consider the coefficient of the monomial
$\prod_{i=1}^r (z_i x_i)^{i-1}$. A simple induction argument shows that in $R(\underline x, \underline z)$, this monomial only appears in the summand indexed by $\sigma = id$, with a coefficient of $(-1)^{r(r-1)/2}.$ In $L(\underline x, \underline z)$ the monomial arises from the product
\[
\prod_{i=1}^r \prod_{l<i} (- \z_l \cdot \x_l ).
\]
We see that the coefficients match, thus the two functions agree.
\end{proof}

\begin{Corollary} For $\delta^-\colon \{1, \dots, r\} \to \{1, \dots, n\}$ an injective function, let $S_r \cdot \delta^-$ denote the set of all possible functions obtained by pre-composing $\delta^-$ with a permutation of $\{1, \dots, r\}.$
The identity
\begin{equation}\label{e:sym2}
\sum_{f^- \in S_r \cdot \delta^-} C^H_{f^-, \delta^+}=C_{\delta^-, \delta^+}
\end{equation}
holds, where $C_{\delta^-, \delta^+}$ was defined in \eqref{CD}.
\end{Corollary}

\begin{proof}
By reindexing the sets $\{\x_1, \dots, \x_n\}$ and $\{\z_1, \dots, \z_n\}$, we may assume without loss of generality that $\delta^- = \delta^+ = \{1, \dots, r\}$. Then in this case, \eqref{e:sym2} simplifies to the equation
\begin{align*}
&\sum_{\sigma \in S_r} \prod_{k=1}^r \Biggl(\prod_{\stackrel{1 \leq j \leq n}{j\neq {\sigma(k)}}}\frac{ \sin( ( \x_{k}- \z_j)/2 {\rm i}) }{\sin(( \z_{\sigma(k)}- \z_j)/2 {\rm i})} \prod_{i<k} \frac{\sin\bigl(\bigl(\z_{\sigma(k)} - \z_{\sigma(i)}\bigr)/2{\rm i}\bigr)}{ \sin((\x_{k} - \x_{i})/2{\rm i})}\Biggr) \nonumber\\
&\qquad{}= \prod_{k=1}^r \prod_{j > r}\frac{ \sin( ( \x_{k}- \z_j)/2 {\rm i}) }{\sin(( \z_{k}- \z_j )/2 {\rm i})}
\end{align*}
or, equivalently,
\begin{align}
&\sum_{\sigma \in S_r} \prod_{k=1}^r \Biggl(\prod_{\stackrel{1 \leq j \leq r}{j\neq {\sigma(k)}}}\frac{ \sin( ( \x_{k}- \z_j)/2 {\rm i}) }{\sin(( \z_{\sigma(k)}- \z_j)/2 {\rm i})} \prod_{i<k} \frac{\sin((\z_{\sigma(k)} - \z_{\sigma(i)})/2{\rm i})}{ \sin((\x_{k} - \x_{i})/2{\rm i})}\Biggr)
= 1 .
\label{e:sym3}
\end{align}
Observe that
\begin{align*}
& \prod_{k=1}^r\Biggl( \prod_{\stackrel{1 \leq j \leq r}{j\neq {\sigma(k)}}} \frac{ 1}{\sin(( \z_{\sigma(k)}- \z_j)/2 {\rm i})}
 \prod_{i<k} \sin((\z_{\sigma(k)} - \z_{\sigma(i)})/2{\rm i}) \Biggr)\\
&\qquad{} = \prod_{k=1}^r \prod_{j>k} \frac{ 1}{\sin(( \z_{\sigma(k)}- \z_{\sigma(j)})/2 {\rm i})}
 = \op{sgn}(\sigma) \prod_{k=1}^r \prod_{j<k} \frac{ 1}{\sin(( \z_{j}- \z_{k})/2 {\rm i})}
\end{align*}
and \eqref{e:sym3} is then equivalent to the equality
\begin{align}
\sum_{\sigma \in S_r} \op{sgn}(\sigma) \prod_{k=1}^r \prod_{\stackrel{1 \leq j \leq r}{j\neq {\sigma(k)}}} \sin( ( \x_{k}- \z_j)/2 {\rm i})
= \prod_{k=1}^r \prod_{i<k} \sin(( \z_{i}- \z_{k})/2 {\rm i}) \sin((\x_{k} - \x_{i})/2{\rm i}).
\label{e:sym4}
\end{align}
Noting that
\[
\sin((a-b)/2{\rm i}) = \frac{ {\rm e}^{a/2 + b/2}\bigl( {\rm e}^a - {\rm e}^b\bigr)}{2{\rm i}},
\]
under the change of variables $X_i = {\rm e}^{\x_i}$, $Z_i = {\rm e}^{\z_i}$, \eqref{e:sym4} becomes
\begin{align*}
\prod_{k=1}^r (Z_k X_k)^{(r-1)/2}\sum_{\sigma \in S_r} \op{sgn}(\sigma) \prod_{\stackrel{1 \leq j \leq r}{j\neq {\sigma(k)}}} (X_k - Z_j)
= \prod_{k=1}^r (Z_k X_k)^{(r-1)/2}\prod_{i<k} (Z_i - Z_k) (X_k - X_i).
\end{align*}
Lemma~\ref{l:sym} then gives the desired equality.
\end{proof}

\subsection{Main theorem}

Combining \eqref{e: CHspecial} with \eqref{e:sym2}, we conclude that
\begin{equation}\label{e: Cextraspecial} \widetilde{H^{+}}_{\beta, \delta^+}|_{q_{1, +} = \cdots =q_{r, +} = q_+} = \sum_{\delta^- \in D_-} C_{\delta^-, \delta^+} H^{-}_{\beta, \delta^-}|_{q_{1, -} = \cdots = q_{r, -} = q_-},\end{equation}
where the path of analytic continuation is along $\hat \gamma$ from Proposition~\ref{p:ancon} and the sum is now over all \emph{unordered} size $r$ subsets of $\{1, \dots, n\}$.

Let $\UU_H\colon H^*_\torus(X_{-})\otimes_{R_\torus} S_\torus \to H^*_\torus(X_{+})\otimes_{R_\torus} S_\torus$ be the map defined in \eqref{e:UH}.
Then \eqref{e: Cextraspecial} implies that
\begin{equation}\label{e:ms}\UU_H H_{X_{-}, \beta}|_{q_{1, -} = \cdots =q_{r, -} = q_-} = \widetilde{H_{X_{+}, \beta}}|_{q_{1, +} = \cdots =q_{r, +} = q_+} \end{equation}
for all $\beta$, where the path of analytic continuation is along $\hat \gamma$ and we set $q_+^{-1} = q_-$.

\begin{Theorem} \label{t:st3}
The linear transformation $\UU$ is symplectic, has a well-defined non-equivariant limit, and satisfies
\[
\UU I_{X_{-}} =\widetilde{I_{X_{+}}},
\]
where $\widetilde{I_{X_{+}}}$ denotes the analytic continuation of $I_{X_{+}}$ along the path $\hat \gamma$ and we set $q_+^{-1} = q_-$.
Furthermore, the following diagram commutes:
\[
\begin{tikzcd}
K^0_\torus(X_{-}) \ar[r, "\FM"] \ar[d, "\Psi_-"] & K^0_\torus(X_{+}) \ar[d, "\Psi_+"] \\
\widetilde{\cc H}_{X_{-}} \ar[r, "\UU"] & \widetilde{\cc H}_{X_{+}}.
\end{tikzcd}
\]
\end{Theorem}
\begin{proof}
Equation \eqref{ItoH3} and the analogous statement for $X_-$ together with \eqref{e:ms} immediately implies that $\UU I_{X_{-}} =\widetilde{I_{X_{+}}}.$ The fact that $\UU$ is symplectic and has a well-defined non-equivariant limit is Corollary~\ref{c:symp}. The commuting diagram the follows from Proposition~\ref{p:Kcom}.
\end{proof}

\begin{Remark}[big $I$-functions and the comparison with \cite{LSW}]\label{r:bigcomp}
In the special case that $B$ is a~point, the wall crossing from $X_+$ to $X_-$ is the standard Grassmann flop.
In this case, by Webb \cite[Section~6.3]{Web}, one can obtain \emph{big} $I$-functions of $X_{\pm}$, denoted $\mathbb I_{X_{\pm}}(q_\pm, \mathbf x, z)$, as a modification of $I_{X_\pm}(q_\pm, z)$.
The respective big $I$-functions generate the entire Lagrangian cones $\cc L_{X_\pm}$, and therefore fully determine the genus-zero Gromov--Witten theory of $X_\pm$. A straightforward generalization of the computations above shows that the transformation $\UU$ from Theorem~\ref{t:st3} also identifies these big $I$-functions:
\begin{equation}\label{e:biggie} \UU \mathbb I_{X_{-}}(q_- , \mathbf x,z) =\widetilde{ \mathbb I_{X_{+}}}(q_+ , \mathbf x,z),\end{equation}
where $\widetilde{ \mathbb I_{X_{+}}}$ denotes the analytic continuation in $q_+$ as in Theorem~\ref{t:st3}.

In the case that $B$ is a point,
Lutz--Shafi--Webb \cite{LSW} have obtained a different proof that there exists a symplectic transformation, which we will denote here by $\UU^{\rm LSW}$, identifying the big $I$-functions of $X_\pm$ after the analytic continuation.
Rather than comparing to the Fourier--Mukai transform $\FM$ however, the symplectic transformation in $\UU^{\rm LSW}$ is obtained directly from wall crossing of the associated abelian GIT quotients $X_{T, \pm}$. Their result is related to Theorem~\ref{t:st3} as follows. Let $\mathbb I_{X_\pm}(q_\pm, \mathbf x ,z)$ be the big $I$-functions of the previous paragraph, and let
$g^* \mathbb I^{\pm}_G(\mathbf y_\pm, \mathbf x ,z)$ denote the $I$-functions used in \cite{LSW}. Comparing explicit formulas, we have
\begin{equation}\label{Icomp} \mathbb I_{X_\pm}(q_\pm , \mathbf x,z)|_{q_\pm = \mathbf y_\pm} = \mathbb I^{\pm}_G(\mathbf y_\pm , \mathbf x,z)|_{\mathbf q = 1}.
\end{equation}

Under the identification
$q_\pm \mapsto \mathbf y_\pm$, the path $\delta$ of analytic continuation used in loc.~cit.\ matches that given in Proposition~\ref{p:ancon}. In \cite{LSW}, the abelian/non-abelian correspondence is used to show that there is a \emph{unique} linear map $\UU^{\rm LSW}$ satisfying
\[
\UU^{\rm LSW} \mathbb I^{-}_G(\mathbf y_- \mathbf x,z)|_{\mathbf q = 1} = \widetilde{\mathbb I^{+}_G}(\mathbf y_+ \mathbf x,z)|_{\mathbf q = 1}.
\]
From this, \eqref{Icomp} and \eqref{e:biggie}, we conclude that $\UU = \UU^{\rm LSW}$, thus providing a consistency check between the two results.
\end{Remark}

\subsection{Central charge}

Motivated by the notion of central charge as given, e.g., in \cite{AL, IriInt}, we make a similar definition. The above wall crossing result takes an especially simple form using this language.
\begin{Definition}

For $E \in K^0_\torus(X_\pm)$, define the \emph{relative quasimap central charge}
\begin{align*}
Z_{X_\pm}(E)(Q, q_\pm, z) := \langle I_{X_{\pm}}(Q, q_\pm, -z), \Psi_\pm(E)\rangle_{X_\pm}.
\end{align*}
\end{Definition}

With this definition, we obtain as an immediate corollary that the central charge of an element~$E$ and of its Fourier--Mukai transform are related by analytic continuation.
\begin{Corollary}
For $E \in K^0_\torus(X_-)$, the analytic continuation of $Z_{X_+}(\FM(E))(Q, q_+, z)$ to ${q_+\! =\! \infty}$ is $Z_{X_-}(E)(Q, q_-, z)$.
\end{Corollary}
\begin{proof}
We use the diagram in Theorem~\ref{t:st} and the fact that $\UU$ preserves the symplectic pairing,
\begin{align*}
\widetilde{Z_{X_+}}(\FM(E))(Q, q_+, z)& =\bigl\langle \widetilde{I_{X_{+}}}(Q, q_+, -z), \Psi_+(\FM(E))\bigr\rangle_{X_+} \\
& = \langle \UU I_{X_{-}}(Q, q_-, -z), \UU \Psi_-(E)\rangle_{X_+}\\
& = \langle I_{X_{-}}(Q, q_-, -z), \Psi_-(E)\rangle_{X_-}
= Z_{X_-}(E)(Q, q_-, z).
\tag*{\qed}
\end{align*}
\renewcommand{\qed}{}
\end{proof}

\appendix
\section{Proof of Proposition~\ref{p:ac1}}\label{appendix}

Using the identity
\begin{equation}\label{e:Gamma}\Gamma(1-z) \Gamma(z) = \frac{\pi}{\sin(\pi z)},\end{equation}
 the denominators of $H_{X_{+}, \beta}$ may be rewritten as
\[
\frac{ \prod_{i=1}^n \sin( (-\z_i -\y)/2 {\rm i}) )}{\pi^n}
\frac{\prod_{i=1}^n (-1)^{\beta \cdot \MB_i + e}\Gamma((-\z_i - \y)/2 \pi {\rm i} + \beta \cdot \MB_i + e)}
{\prod_{j=1}^n \Gamma(1 + (- \x_j - \y)/2 \pi {\rm i} + \beta \cdot \LB_j+e) }.
\]

Let $B_l^{+}$ denote the $l$th fixed locus, isomorphic to $B$. By Corollary~\ref{c:rest}, the restriction of the relative hyperplane class $-y$ to $B_l^{+}$ is $ \x_l$. Thus $H_{\beta,l}^{+} :=H_{X_{+}, \beta}|_{B_l^{+}}$ is given by
\begin{align}
&\sum_{e \geq -\beta \cdot \LB_l} q_+^{e+ \x_l/2 \pi {\rm i}}\nonumber\\
&\hphantom{\sum_{e \geq -\beta \cdot \LB_l}}{}\times\frac{ \prod_{i=1}^n \sin( (-\z_i+\x_l)/2 {\rm i}) }{\pi^n}\frac{\prod_{i=1}^n (-1)^{\beta \cdot \MB_i + e}\Gamma((-\z_i+\x_l)/2 \pi {\rm i} + \beta \cdot \MB_i + e)}
{\prod_{j=1}^n \Gamma(1 + (- \x_j+\x_l)/2 \pi {\rm i} + \beta \cdot \LB_j+e) }.\label{e:Hp}
\end{align}
This satisfies the differential equation
\[
\left[\prod_{j=1}^n (\partial_q -x_j/2 \pi {\rm i} + \beta \cdot L_j) - q_+(-1)^n\prod_{j=1}^n (\partial_q -z_j/2 \pi {\rm i} + \beta \cdot M_j) \right] \Phi(q_+) = 0,
\]
where $\partial_q = q_{+} \frac{\partial}{\partial q_{+}}$. This differential equation has singularities at $0, \infty,$ and $(-1)^n$. Following \cite{CIJ}, our path $\gamma$ of analytic continuation in the $\log(q_{+})$ will go through $\mathfrak{Im}(\log(q_{+})) = (n-1)\pi {\rm i}$ when $\mathfrak{Re}(\log(q_{+}))$ is close to 0. See Figure~\ref{f1} for a picture of the path.
We will use the Mellin--Barnes method to compute this analytic continuation of $H_{\beta,l}^{+}$ to $q_+ = \infty$.

We may rewrite \eqref{e:Hp} as a contour integral
\begin{align} \label{e:cint}&\frac{1}{2 \pi {\rm i}}\frac{ \prod_{i=1}^n \sin( (-\z_i+\x_l)/2 {\rm i}) }{\pi^n} \\
&\times\! \int_C \Gamma( s) \Gamma(1 - s) q_+^{s + \x_l/2 \pi {\rm i}} {\rm e}^{ -\pi {\rm i} (\sum_{i=1}^n \beta \cdot \MB_i+(n-1)\cdot s)} \frac{ \prod_{i=1}^n \Gamma((-\z_i+\x_l)/2 \pi {\rm i} + \beta \cdot \MB_i + s)}
{\prod_{j=1}^n \Gamma(1 + (- \x_j+\x_l)/2 \pi {\rm i} + \beta \cdot \LB_j+s) },
\nonumber \end{align}
where $C$ encloses all integers $e \geq -\beta \cdot \LB_l$.

Closing the curve on the other side, we obtain residues at integers $d < -\beta \cdot \LB_l$ and at
\[
- (-\z_k+\x_l)/2 \pi {\rm i} - d,\qquad 1 \leq k \leq n, d \in \ZZ_{\geq \beta \cdot \MB_k}.
\]
The residues at integers $d \leq -\beta \cdot \LB_l$ vanish. The integral is then equal to minus the sum of the remaining residues:
\begin{gather*} 
 -\sum_{1 \leq k \leq n}
\sum_{d \geq \beta \cdot \MB_k}\Gamma(- (-\z_k+\x_l)/2 \pi {\rm i} -d) \Gamma(1 + (-\z_k+\x_l)/2 \pi {\rm i} + d) q_+^{ \z_k/2 \pi {\rm i}-d } \\ \nonumber
\times
\frac{ {\rm e}^{- \pi {\rm i} ( \sum_{i=1}^n \beta \cdot \MB_i + (n-1)\cdot ( - (-\z_k+\x_l)/2 \pi {\rm i} - d)+d- \beta \cdot \MB_k)}}{(d - \beta \cdot \MB_k)!}\frac{\prod_{i\neq k} \Gamma((-\z_i + \z_k)/2 \pi {\rm i} + \beta \cdot \MB_i -d)}
{\prod_{j=1}^n \Gamma(1 + (- \x_j+\z_k)/2 \pi {\rm i} + \beta \cdot \LB_j -d) }\\ \nonumber
\times
 \frac{ \prod_{i=1}^n \sin( (-\z_i+\x_l)/2 {\rm i}) }{\pi^n}
 \\ \nonumber
= -\sum_{1 \leq k \leq n}
\sum_{d \geq \beta \cdot \MB_k}\Gamma(- (-\z_k+\x_l)/2 \pi {\rm i} ) \Gamma(1 + (-\z_k+\x_l)/2 \pi {\rm i}) q_+^{ \z_k/2 \pi {\rm i}-d } \\ \nonumber
\times
\frac{ {\rm e}^{ (n-1) ( \x_l-\z_k)/2}}
{(d - \beta \cdot \MB_k)!\prod_{j=1}^n \Gamma(1 + (- \x_j+\z_k)/2 \pi {\rm i} + \beta \cdot \LB_j -d) }\\ \nonumber
\times
 \frac{ \prod_{i=1}^n \sin( (-\z_i+\x_l)/2 {\rm i}) }{\prod_{i \neq k} \Gamma(1 - (-\z_i +\z_k)/2 \pi {\rm i} - \beta \cdot \MB_i + d) \cdot \pi \prod_{i \neq k}\sin((-\z_i + \z_k)/2 {\rm i})}
\\ \nonumber
 = \sum_{1 \leq k \leq n}
\sum_{d \geq \beta \cdot \MB_k} {\rm e}^{ (n-1) ( \x_l-\z_k)/2}\frac{ \prod_{i\neq k} \sin( (-\z_i+\x_l)/2 {\rm i}) }{\prod_{i \neq k}\sin((-\z_i + \z_k)/2 {\rm i})} \\
\times
\frac{q_+^{ \z_k/2 \pi {\rm i}-d } }
{\prod_{j=1}^n \Gamma(1 + (- \x_j+\z_k)/2 \pi {\rm i} + \beta \cdot \LB_j -d) \prod_{i =1}^n \Gamma(1 - (-\z_i +\z_k)/2 \pi {\rm i} - \beta \cdot \MB_i + d)}, \nonumber
\end{gather*}
where we have again simplified using the Gamma function identity \eqref{e:Gamma}.

Let $B_k^{-}$ denote the $k$th fixed locus of $X_{-}$. By \eqref{e:HXm}, observe
that $H_{\beta, k}^{-} := H_{X_{-}, \beta}|_{B_k^{-}}$ is given by
\[
\sum_{d \geq \beta \cdot \MB_k} \frac{ q_-^{-\z_k/2 \pi {\rm i}+d}}
{\prod_{j=1}^n \Gamma(1 + (- \x_j +\z_k)/2 \pi {\rm i} + \beta \cdot \LB_j-d)\Gamma(1 + (\z_j -\z_k)/2 \pi {\rm i} - \beta \cdot \MB_j + d)}.
\]
The result then follows immediately from \eqref{eU1}.

\subsection*{Acknowledgements}
We would like to thank Tom Coates, Wendelin Lutz, Jeongseok Oh, Yongbin Ruan, Ed Segal, Qaasim Shafi, and Rachel Webb for the useful discussions. We also thank the anonymous referees for valuable comments and suggestions on a previous draft.
N.~Priddis is supported by the Simons Foundation Travel Grant 586691. M.~Shoemaker is supported by the Simons Foundation Travel Grant 958189.
Y.~Wen is supported by a KIAS Individual Grant (MG083902) at the Korea Institute for Advanced Study.


\begin{thebibliography}{99}
\footnotesize\itemsep=0pt

\bibitem{AL}
Aleshkin K., Liu C.-C.M., {H}iggs--{C}oulomb correspondence and wall-crossing in
 abelian {GLSM}s, \href{https://arxiv.org/abs/2301.01266}{arXiv:2301.01266}.

\bibitem{AB}
Atiyah M.F., Bott R., The moment map and equivariant cohomology,
 \href{https://doi.org/10.1016/0040-9383(84)90021-1}{\textit{Topology}} \textbf{23} (1984), 1--28.

\bibitem{BCFMV}
Ballard M.R., Chidambaram N.K., Favero D., McFaddin P.K., Vandermolen R.R.,
 Kernels for {G}rassmann flops, \href{https://doi.org/10.1016/j.matpur.2021.01.005}{\textit{J.~Math. Pures Appl.}} \textbf{147}
 (2021), 29--59, \href{https://arxiv.org/abs/1904.12195}{arXiv:1904.12195}.

\bibitem{BCK}
Bertram A., Ciocan-Fontanine I., Kim B., Two proofs of a conjecture of {H}ori
 and {V}afa, \href{https://doi.org/10.1215/S0012-7094-04-12613-2}{\textit{Duke Math.~J.}} \textbf{126} (2005), 101--136,
 \href{https://arxiv.org/abs/math.AG/0304403}{arXiv:math.AG/0304403}.

\bibitem{BCK2}
Bertram A., Ciocan-Fontanine I., Kim B., Gromov--{W}itten invariants for
 abelian and nonabelian quotients, \href{https://doi.org/10.1090/S1056-3911-07-00456-0}{\textit{J.~Algebraic Geom.}} \textbf{17}
 (2008), 275--294, \href{https://arxiv.org/abs/math.AG/0407254}{arXiv:math.AG/0407254}.

\bibitem{Brion}
Brion M., Equivariant cohomology and equivariant intersection theory, in
 Representation Theories and Algebraic Geometry ({M}ontreal, {PQ}, 1997),
 \textit{NATO Adv. Sci. Inst. Ser.~C: Math. Phys. Sci.}, Vol. 514, \href{https://doi.org/10.1007/978-94-015-9131-7_1}{Kluwer},
 Dordrecht, 1998, 1--37, \href{https://arxiv.org/abs/math.AG/9802063}{arXiv:math.AG/9802063}.

\bibitem{Brown}
Brown J., Gromov--{W}itten invariants of toric fibrations, \href{https://doi.org/10.1093/imrn/rnt030}{\textit{Int. Math.
 Res. Not.}} \textbf{2014} (2014), 5437--5482, \href{https://arxiv.org/abs/0901.1290}{arXiv:0901.1290}.

\bibitem{BG}
Bryan J., Graber T., The crepant resolution conjecture, in Algebraic
 Geometry---{S}eattle 2005. {P}art~1, \textit{Proc. Sympos. Pure Math.},
 Vol.~80, \href{https://doi.org/10.1090/pspum/080.1/2483931}{American Mathematical Society}, Providence, RI, 2009, 23--42,
 \href{https://arxiv.org/abs/math.AG/0610129}{arXiv:math.AG/0610129}.

\bibitem{BLV}
Buchweitz R.-O., Leuschke G.J., Van~den Bergh M., Non-commutative
 desingularization of determinantal varieties, {II}: {A}rbitrary minors,
 \href{https://doi.org/10.1093/imrn/rnv207}{\textit{Int. Math. Res. Not.}} \textbf{2016} (2016), 2748--2812,
 \href{https://arxiv.org/abs/1106.1833}{arXiv:1106.1833}.

\bibitem{ChIR}
Chiodo A., Iritani H., Ruan Y., Landau--{G}inzburg/{C}alabi--{Y}au
 correspondence, global mirror symmetry and {O}rlov equivalence, \href{https://doi.org/10.1007/s10240-013-0056-z}{\textit{Publ.
 Math. Inst. Hautes \'Etudes Sci.}} \textbf{119} (2014), 127--216,
 \href{https://arxiv.org/abs/1201.0813}{arXiv:1201.0813}.

\bibitem{ChR}
Chiodo A., Ruan Y., Landau--{G}inzburg/{C}alabi--{Y}au correspondence for
 quintic three-folds via symplectic transformations, \href{https://doi.org/10.1007/s00222-010-0260-0}{\textit{Invent. Math.}}
 \textbf{182} (2010), 117--165, \href{https://arxiv.org/abs/0812.4660}{arXiv:0812.4660}.

\bibitem{CKM14}
Ciocan-Fontanine I., Kim B., Maulik D., Stable quasimaps to {GIT} quotients,
 \href{https://doi.org/10.1016/j.geomphys.2013.08.019}{\textit{J.~Geom. Phys.}} \textbf{75} (2014), 17--47, \href{https://arxiv.org/abs/1106.3724}{arXiv:1106.3724}.

\bibitem{CKS}
Ciocan-Fontanine I., Kim B., Sabbah C., The abelian/nonabelian correspondence
 and {F}robenius manifolds, \href{https://doi.org/10.1007/s00222-007-0082-x}{\textit{Invent. Math.}} \textbf{171} (2008),
 301--343, \href{https://arxiv.org/abs/math.AG/0610265}{arXiv:math.AG/0610265}.

\bibitem{CIJ}
Coates T., Iritani H., Jiang Y., The crepant transformation conjecture for
 toric complete intersections, \href{https://doi.org/10.1016/j.aim.2017.11.017}{\textit{Adv. Math.}} \textbf{329} (2018),
 1002--1087, \href{https://arxiv.org/abs/1410.0024}{arXiv:1410.0024}.

\bibitem{CIJS}
Coates T., Iritani H., Jiang Y., Segal E., {$K$}-theoretic and categorical
 properties of toric {D}eligne--{M}umford stacks, \href{https://doi.org/10.4310/PAMQ.2015.v11.n2.a3}{\textit{Pure Appl. Math.~Q.}}
 \textbf{11} (2015), 239--266, \href{https://arxiv.org/abs/1410.0027}{arXiv:1410.0027}.

\bibitem{CIT}
Coates T., Iritani H., Tseng H.-H., Wall-crossings in toric {G}romov--{W}itten
 theory.~{I}. {C}repant examples, \href{https://doi.org/10.2140/gt.2009.13.2675}{\textit{Geom. Topol.}} \textbf{13} (2009),
 2675--2744, \href{https://arxiv.org/abs/math.AG/0611550}{arXiv:math.AG/0611550}.

\bibitem{CLS}
Coates T., Lutz W., Shafi Q., The abelian/nonabelian correspondence and
 {G}romov--{W}itten invariants of blow-ups, \href{https://doi.org/10.1017/fms.2022.46}{\textit{Forum Math. Sigma}}
 \textbf{10} (2022), e67, 33~pages, \href{https://arxiv.org/abs/2108.10922}{arXiv:2108.10922}.

\bibitem{CR}
Coates T., Ruan Y., Quantum cohomology and crepant resolutions: a~conjecture,
 \href{https://doi.org/10.5802/aif.2766}{\textit{Ann. Inst. Fourier (Grenoble)}} \textbf{63} (2013), 431--478,
 \href{https://arxiv.org/abs/0710.5901}{arXiv:0710.5901}.

\bibitem{EG}
Edidin D., Graham W., Riemann--{R}och for equivariant {C}how groups,
 \href{https://doi.org/10.1215/S0012-7094-00-10239-6}{\textit{Duke Math.~J.}} \textbf{102} (2000), 567--594,
 \href{https://arxiv.org/abs/math.AG/9905081}{arXiv:math.AG/9905081}.

\bibitem{Fulton}
Fulton W., Intersection theory, \textit{Ergeb. Math. Grenzgeb.~(3)}, Vol.~2,
 \href{https://doi.org/10.1007/978-3-662-02421-8}{Springer}, Berlin, 1984.

\bibitem{GW}
Gonz\'alez E., Woodward C.T., A wall-crossing formula for {G}romov--{W}itten
 invariants under variation of git quotient, \href{https://doi.org/10.1007/s00208-023-02622-w}{\textit{Math. Ann.}} \textbf{388}
 (2024), 4135--4199, \href{https://arxiv.org/abs/1208.1727}{arXiv:1208.1727}.

\bibitem{IriInt}
Iritani H., An integral structure in quantum cohomology and mirror symmetry for
 toric orbifolds, \href{https://doi.org/10.1016/j.aim.2009.05.016}{\textit{Adv. Math.}} \textbf{222} (2009), 1016--1079,
 \href{https://arxiv.org/abs/0903.1463}{arXiv:0903.1463}.

\bibitem{ILLW}
Iwao Y., Lee Y.-P., Lin H.-W., Wang C.-L., Invariance of {G}romov--{W}itten theory
 under a~simple flop, \href{https://doi.org/10.1515/CRELLE.2011.097}{\textit{J.~Reine Angew. Math.}} \textbf{663} (2012),
 67--90, \href{https://arxiv.org/abs/0804.3816}{arXiv:0804.3816}.

\bibitem{KKP}
Katzarkov L., Kontsevich M., Pantev T., Hodge theoretic aspects of mirror
 symmetry, in From {H}odge Theory to Integrability and {TQFT} tt*-geometry,
 \textit{Proc. Sympos. Pure Math.}, Vol.~78, \href{https://doi.org/10.1090/pspum/078/2483750}{American Mathematical Society},
 Providence, RI, 2008, 87--174, \href{https://arxiv.org/abs/0806.0107}{arXiv:0806.0107}.

\bibitem{LLW}
Lee Y.-P., Lin H.-W., Wang C.-L., Flops, motives, and invariance of quantum rings,
 \href{https://doi.org/10.4007/annals.2010.172.243}{\textit{Ann. of Math.}} \textbf{172} (2010), 243--290,
 \href{https://arxiv.org/abs/math.AG/0608370}{arXiv:math.AG/0608370}.

\bibitem{LLW2}
Lee Y.-P., Lin H.-W., Wang C.-L., Invariance of quantum rings under ordinary
 flops~{II}: {A}~quantum {L}eray--{H}irsch theorem, \href{https://doi.org/10.14231/AG-2016-027}{\textit{Algebr. Geom.}}
 \textbf{3} (2016), 615--653, \href{https://arxiv.org/abs/1311.5725}{arXiv:1311.5725}.

\bibitem{LR}
Li A.-M., Ruan Y., Symplectic surgery and {G}romov--{W}itten invariants of
 {C}alabi--{Y}au 3-folds, \href{https://doi.org/10.1007/s002220100146}{\textit{Invent. Math.}} \textbf{145} (2001),
 151--218, \href{https://arxiv.org/abs/math.AG/9803036}{arXiv:math.AG/9803036}.

\bibitem{LSW}
Lutz W., Shafi Q., Webb R., Crepant transformation conjecture for the
 {G}rassmannian flop, \href{https://arxiv.org/abs/2404.12302}{arXiv:2404.12302}.

\bibitem{Oh1}
Oh J., Quasimaps to {GIT} fibre bundles and applications, \href{https://doi.org/10.1017/fms.2021.48}{\textit{Forum Math.
 Sigma}} \textbf{9} (2021), e56, 39~pages, \href{https://arxiv.org/abs/1607.08326}{arXiv:1607.08326}.

\bibitem{Oko15}
Okounkov A., Lectures on {K}-theoretic computations in enumerative geometry, in
 Geometry of Moduli Spaces and Representation Theory, \textit{IAS/Park City
 Math. Ser.}, Vol.~24, \href{https://doi.org/10.1090/pcms/024}{American Mathematical Society}, Providence, RI, 2017,
 251--380, \href{https://arxiv.org/abs/1512.07363}{arXiv:1512.07363}.

\bibitem{Ru}
Ruan Y., Surgery, quantum cohomology and birational geometry, in Northern
 {C}alifornia {S}ymplectic {G}eometry {S}eminar, \textit{Amer. Math. Soc.
 Transl. Ser.~2}, Vol. 196, \href{https://doi.org/10.1090/trans2/196/09}{American Mathematical Society}, Providence, RI,
 1999, 183--198, \href{https://arxiv.org/abs/math.AG/9810039}{arXiv:math.AG/9810039}.

\bibitem{Thomason}
Thomason R.W., Une formule de {L}efschetz en {$K$}-th\'eorie \'equivariante
 alg\'ebrique, \href{https://doi.org/10.1215/S0012-7094-92-06817-7}{\textit{Duke Math.~J.}} \textbf{68} (1992), 447--462.

\bibitem{Web}
Webb R., The abelian-nonabelian correspondence for {$I$}-functions,
 \href{https://doi.org/10.1093/imrn/rnab305}{\textit{Int. Math. Res. Not.}} \textbf{2023} (2023), 2592--2648,
 \href{https://arxiv.org/abs/1804.07786}{arXiv:1804.07786}.

\bibitem{Web2}
Webb R., Abelianization and quantum {L}efschetz for orbifold quasimap
 {$I$}-functions, \href{https://doi.org/10.1016/j.aim.2024.109489}{\textit{Adv. Math.}} \textbf{439} (2024), 109489, 59~pages,
 \href{https://arxiv.org/abs/2109.12223}{arXiv:2109.12223}.

\end{thebibliography}

\pdfbookmark[1]{References}{ref}
\LastPageEnding

\end{document}